\def\year{2022}\relax
\documentclass[letterpaper]{article} 
\usepackage{aaai22}  
\usepackage{times}  
\usepackage{helvet}  
\usepackage{courier}  
\usepackage[hyphens]{url}  
\usepackage{graphicx} 
\usepackage{natbib}  
\usepackage{caption} 
\DeclareCaptionStyle{ruled}{labelfont=normalfont,labelsep=colon,strut=off} 
\usepackage{algorithm}
\usepackage{algorithmic}
\usepackage{booktabs}

\usepackage{color}
\usepackage{bbm}
\usepackage{amsmath}
\usepackage{subcaption}
\usepackage{multirow}
\usepackage{lipsum}
\usepackage{xcolor}
\usepackage{sidecap}

\newcommand{\sun}[1]{{\color{black} #1}}

\newcommand{\shi}[1]{{\color{black} #1}}

\DeclareMathOperator*{\argmax}{arg\,max}

\newtheorem{theorem}{Theorem}

\newtheorem{proposition}{Proposition}

\usepackage{newfloat}
\usepackage{listings}
\floatstyle{ruled}
\newfloat{listing}{tb}{lst}{}
\floatname{listing}{Listing}
\title{\sun{Constrained Prescriptive Trees via Column Generation}}
\author {
    Shivaram Subramanian\equalcontrib,
    Wei Sun\equalcontrib,
    Youssef Drissi,
    Markus Ettl
}
\affiliations {
    IBM Research, Yorktown Heights, NY 10591, USA\\
    subshiva, sunw, youssefd, mettl@us.ibm.com
}

\usepackage{bibentry}

\begin{document}

\maketitle

\begin{abstract}
With the abundance of available data, many enterprises seek to implement data-driven prescriptive analytics to help them make informed decisions. These prescriptive policies need to satisfy operational constraints, and proactively eliminate rule conflicts, both of which are ubiquitous in practice. It is also desirable for them to be simple and interpretable, so they can be easily verified and implemented. Existing approaches from the literature center around constructing variants of prescriptive decision trees to generate interpretable policies. However, none of the existing methods are able to handle constraints.  In this paper, we propose a scalable method that solves the constrained prescriptive policy generation problem. We introduce a  novel path-based mixed-integer program (MIP) formulation which identifies \shi{a (near)} optimal policy efficiently via column generation. The policy generated can be represented as a multiway-split tree which \shi{is} more interpretable and informative than \shi{a binary-split tree} due to its shorter rules. We demonstrate the efficacy of our method with extensive  experiments on both synthetic and real datasets.\end{abstract}

\section*{Introduction}
With the abundance of available data, many enterprises seek to implement data-driven prescriptive analytics to help them make better decisions. \sun{Unlike \emph{predictive} analytics, which applies mathematical models to forecast potential  outcomes, \emph{prescriptive} analytics \shi{determines} the best course of action for the future, frequently utilizing the outcome of a predictive model  (\citealp{lepenioti2020prescriptive})}. 
For instance, in revenue management,  pricing decisions are based on estimated demand learnt from historical sales data; 
in healthcare, personalized medicine is prescribed based on past patients' responses towards different treatment options. A key characteristic of prescriptive analytics is the presence of constraints, which is largely missing from predictive analytics. Constraints are  ubiquitous in practice: 
a revenue-maximization problem may limit prices such that the resulting demand does not exceed production capacity; similarly, a patient's available treatment options may be limited by his or her pre-existing conditions.

Despite machine learning (ML) making huge strides in recent years, there are several obstacles standing in the way of widespread adoption of prescriptive analytics in practice. Firstly, prescriptive analytics \sun{often} works in conjunction with predictive models, whose output may be nonlinear, non-convex and/or discontinuous (and therefore not differentiable), 
\shi{resulting in intractable downstream optimization problems.}
To generate meaningful policies, 
it is critical to embed the myriad of business requirements  into the \shi {decision} optimization problem, exacerbating the computational challenge at hand.   Secondly, powerful black-box predictive models inevitably obscure the downstream decision-making process, making it difficult  for enterprises to understand and trust the prescribed policies. Meanwhile, recommendations based on these models tend to be highly complex, making implementation and  maintenance of these policies challenging and cumbersome.

There have been some  recent works studying the problem of generating interpretable prescriptive policy. 
A common theme is to construct variants of decision trees which represent the prescribed policies, as trees provide  ``human-friendly'' explanations (\citealp{frosst2017distilling,miller2019explanation}). Each path in a tree from the root node to a leaf node  corresponds to a policy, and all samples in a leaf are prescribed with a same action.  
The procedure of generating such a tree-based policy may vary - the prediction step can be either embedded in the policy generation step (\citealp{kallus2017recursive,bertsimas2019optimal}), or explicitly separated from  prescription (\citealp{zhou2018offline,amram2020optimal,biggs2020model}). Meanwhile, the prescriptive tree can be constructed either greedily (\citealp{zhou2018offline,biggs2020model}) or optimally (\citealp{kallus2017recursive,amram2020optimal}). However, none of the existing methods \shi{are} capable of handling operational constraints which are paramount for successful enterprise adoption of prescriptive analytics - this is the gap which we intend to fill with this work. 

In this paper, we propose a scalable method to solve the constrained prescriptive policy generation problem. We consider the \emph{predictive teacher} and \emph{prescriptive student} framework  \cite{zhou2018offline,amram2020optimal,biggs2020model}, where a predictive model (teacher) is first trained to produce counterfactual outcomes associated with different actions for every sample in the dataset. Based on the counterfactual estimations, the downstream prescriptive model (student) determines the best set of policies to divide the samples that optimizes the  given objective. 
Our key contributions are as follows:

\textbf{1. Constrained policy prescription} To the best of our knowledge, we are the first to propose a prescriptive tree which is capable of handling  constraints that span across multiple rules. Existing tree-based approaches in the literature typically produce unconstrained policies that ignore rule conflicts and operational feasibility, severely limiting their usefulness in practice.

\textbf{2. Novel MIP formulation}
We model rules or paths in a tree explicitly as it provides a natural way to impose  constraints unlike prior MIP methods which rely on arc-based feature-level models. 
More specifically, we show that the constrained tree-based policy prescription 
problem can be formulated as a set-partitioning problem with side-constraints. A 
key challenge here is that the search space for rules may be computationally 
prohibitive to solve the problem directly for high dimensional datasets.  To 
overcome this hurdle, we implement column generation \cite{ford1958suggested}, where promising candidate rules can be generated dynamically. Column 
Generation (CG) can be applied as an exact method as well as an efficient 
technique to obtain near-optimal solutions in practical settings.

\textbf{3. Scalable algorithm}
The aforementioned \emph{path-based} MIP formulation enables breakthroughs in 
scalability over prior \emph{arc-based} MIP decision tree methods (e.g., \citealp{bertsimas2019optimal,amram2020optimal}), which suffer from a scalability issue 
as their formulation requires $O(2^kM)$ binary variables where $k$ and $M$ refer to the depth of the tree and the number of samples respectively. In contrast, the 
number of binary decision variables in our MIP formulation is independent of the 
sample size and equals the number of candidate decision rules. While this number 
can also be prohibitive in the worst case, we only need to generate candidates as 
needed to find a near-optimal solution and typically the resultant number of 
candidate rules is far less than the number of possible rules. As a result, while 
existing approaches struggle to process datasets beyond $10^5$ samples and depth 
2, a first crude implementation of our proposed CG heuristic algorithm handles 
$10^6$ samples at a depth equivalent of 5, with the potential to solve larger datasets 
using more efficient implementations.

\textbf{4. More interpretable policies} Contrary to the existing tree-based methods which typically construct a binary-split tree where each node has at most two children, we consider a multiway-split tree, where a node may have more than two \shi{child} nodes (refer to Figure~\ref{multiway_split} for an example). Multiway trees offer the advantage over binary trees that an attribute rarely appears more than once in any path from root to leaf, which are easier to comprehend than its binary counterparts (\citealp{fulton1995efficient}).

\textbf{5. Extensive numerical results} We perform extensive experiments to show that the proposed approach produces solutions of superior quality compared to prior prescriptive tree approaches in the literature. We present several use cases with public datasets to demonstrate our method's flexibility in handling a variety of constraints.

\section*{Related Literature}

There has been a surge of interest in making ML models more interpretable  in the recent years (e.g., \citealp{ribeiro2016should,lundberg2017unified,guidotti2018survey,carvalho2019machine}). 
A common framework used to gain interpretability is \emph{knowledge distillation} (\citealp{hinton2015distilling}), where a complex model serves as a teacher and a \shi{simpler} model learns to mimic the teacher as a student. For instance, this concept has  been implemented in explainable reinforcement learning (XRL) \cite{liu2018toward,puiutta2020explainable}, where a  black-box teacher model first learns the policies, then a regression tree is trained to approximate these policies. A notable distinction of the  framework considered in this work lies in the different roles played by the teacher and the student models. More specifically, in our setting, the prescriptive student does not 
 merely replicate the teacher's behavior, rather it determines the optimal policies. Meanwhile, the teacher  in our framework does not produce policy, rather it provides the counterfactual estimations which guide the student in determining the optimal policies.

 As this framework 
explicitly decouples prediction from prescription, one can take advantage of the latest advances in ML and utilize powerful models such as boosted trees or deep neural networks for the prediction task. It is in  contrast to an alternative approach, i.e., embedding prediction inside prescription by combining the tasks of estimating a predictive model and learning the optimal policy  (\citealp{kallus2017recursive,bertsimas2019optimal}). While this approach offers the benefit of learning prediction and prescription from data in a single step, it also limits the predictive models to  piecewise-constant or piecewise-linear functions for computational tractability. Such restriction may lead to severe model misspecification when the underlying structure of the outcome takes a more complex form.

Besides the aforementioned interpretable prescription methods, there exist  approaches from causal analysis, where the focus is  on estimating heterogeneous treatment effects which can in turn quantify the impact of prescribing actions (e.g., \citealp{shalit2017estimating,wager2018estimation,kunzel2019metalearners}). These methods are typically limited to binary actions and are less concerned with interpretability.


Decision trees have been used extensively in machine learning, mainly due to their transparency which allows users to derive interpretable results (\citealp{NEURIPS2020_1373b284}). As learning an optimal  decision tree is NP-complete (\citealp{laurent1976constructing}), popular algorithms such as CART (\citealp{breiman1984classification}) have relied on greedy heuristics to construct trees. Recent advances in modern optimization has facilitated a nascent stream of research that leverages mixed-integer programming to train globally-optimal
trees  (e.g., \citealp{bertsimas2017optimal,bertsimas2019optimal,aghaei2019learning,aghaei2020learning,NEURIPS2020_1373b284}). In particular, \citealp{bertsimas2019optimal} and \citealp{amram2020optimal} have applied this MIP  approach to learn a prescription policy. Notable differences between our work and the aforementioned MIP approach  are 1) instead of constructing a binary-split tree, we construct a multiway-split tree which is more interpretable 
(\citealp{fulton1995efficient}). 
2) Instead of an arc-based formulation, we present a path-based MIP formulation to construct the tree, which allows us to model constraints naturally and solve the problem efficiently via column generation (CG). CG has been successfully adopted by the AI community (e.g., \citealp{bach2008exploring,jawanpuria2011efficient}) and in practice to solve  complex discrete optimization models including airline scheduling (\citealp{klabjan2002airline, subramanian2008effective, bront2009column}), vehicle routing (\citealp{chen2006dynamic}), and inventory planning for supply chain, among others (\citealp{desaulniers2006column,xu2019solving}).

\section*{Problem Formulation}
\subsection{Predictive Teacher with Prescriptive Student Tree}

We assume there are $M$ observational data samples $\{(x_i,\pi_i,y_i)\}_{i=1}^M$, where $x_i \in \mathcal{X}^d$ are features,  $ \pi_i$ refers to the action chosen from a discrete set $\Pi$,
and $y_i$ is the uncertain quantity of interest, which can be a discrete label for classification or a continuous quantity for regression.  In a use case of pricing studied in Section~\ref{sect_real}, $x_i$ represents customer features, $\pi_i$ is the price of a product and $y_i \in \{0,1\}$ indicates whether the product was sold $(1)$ or not $(0)$. In another use case on housing upgrade, $x_i$ represents attributes associated with a house, $\pi_i$ is the construction grade which a builder can choose from, and $y_i$ indicates the sale price of a house.


As a standard practice in causal inference literature, we make the ignorability assumption (\citealp{hirano2004propensity}), i.e., there were no unmeasured confounding variables that affected both the choice of decision and the outcome.

We assume an underlying function $f:\mathcal{X}^d \times \Pi$ which maps the  features and an action to the outcome $y$. An optimal policy selects an action $\tau^*(x)$  for each sample based on its features  to maximize the objective function:  
\begin{equation}
\label{true_rev}
\tau^*(x)= \argmax_{\pi \in \Pi} g(f(x,\pi) )
\end{equation}
In pricing, $f(x,\pi)$ represents the probability of purchase given price $\pi$ and customer feature $x$, while  $g=\pi f(x,\pi)$ denotes the expected revenue. 

In reality, the true response function $f(x,\pi)$ is unknown, but can be estimated. We refer to the proxy response function $\hat{f}$ as the teacher model. There are many possible choices for $\hat{f}$. \citealp{biggs2020model} learn a teacher model by solving an empirical risk minimization problem. 
\citealp{zhou2018offline,amram2020optimal}  use doubly-robust estimators  to estimate the counterfactual outcomes 
(\citealp{dudik2011doubly}). The advantage of separating the prediction task from prescription is that one has the flexibility to choose the ``best'' counterfactual inference model available as the teacher for a given application.    

If we substitute the proxy response function $\hat{f}$ into the optimization problem in (\ref{true_rev}), the prescription problem is reduced to enumerate through $\Pi$ to determine the optimal policy. More concretely,  we use  $g_{i,\pi}$ to denote the estimated counterfactual outcome by the teacher model for a given action $\pi$, i.e., $g_{i,\pi} = g(\hat{f}(x_i,\pi))$, and the policy prescription problem becomes $$\tau(x)= \argmax_{\pi \in \Pi}  \sum_{i=1}^M g_{i,\pi}$$Note that since $g_{i,\pi}$ is just an input
to the downstream optimization problem, the task of policy prescription  becomes much simpler, i.e., for each sample $i$, enumerate through $g_{i,\pi}$ and select the best action.  
One caveat of this approach is that the resulting policy may be too complex when the action space is huge, and may not offer any insight into how the prescription is made. 

\begin{figure}
\begin{subfigure}{.5\textwidth}
  \includegraphics[width=0.8\linewidth]{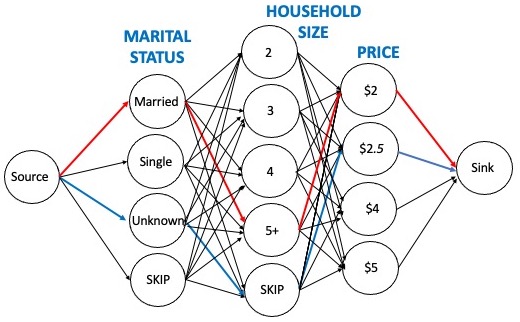}
  \caption{Feature graph with 3  features including action (price)}
  \label{feature_space}
  \end{subfigure}
\begin{subfigure}{.5\textwidth}
  \includegraphics[width=0.8\linewidth]{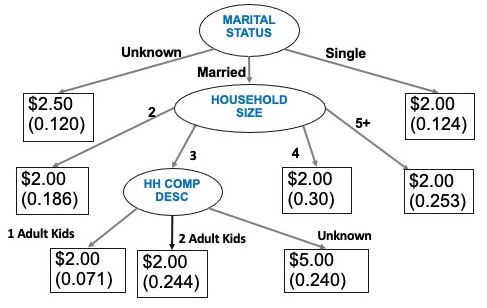}
  \caption{A multiway-split tree with $n=8$ rules}
  \label{multiway_split}
  \end{subfigure}
  \caption{A sample output from the grocery pricing use case }
\end{figure}

To aid interpretability, we restrict the policy to a certain class of prescriptive models. A common choice in the literature is a binary-split decision tree of a pre-specified depth $k$, $\mathcal{T}_b(k)$, which has at most $2^k$ leaves (e.g., \citealp{biggs2020model,amram2020optimal}). In this work, we consider a multiway-split  tree (MT) as the prescriptive student model, denoted as $\mathcal{T}_m(n)$, where $n$ indicates the number of leaves.
An example of a MT with $n=8$ is shown in Fig~\ref{multiway_split}. A MT with $n$ leaves $\mathcal{T}_m(n)$ is comparable to a  binary-split tree grown to its full width with the same number of leaves, i.e.,  $\mathcal{T}_b(k)$ where $2^k = n$. We refer to our model  as \emph{student prescriptive multiway-split tree} problem, or SPMT,
\begin{align}
\textbf{(SPMT)}\qquad \max_{\tau(x) \in \mathcal{T}_m(n)} \sum_{i=1}^M g_{i,\pi} \label{SPMT}
\end{align}
There are many real-life situations where a decision maker may only want to use a subset of features $\hat{x}$ to define a policy. For instance, one may not want to include personal features due to privacy and fairness concerns. This highlights another advantage of the \emph{prescriptive knowledge distillation} framework - an ``informative'' teacher which is trained on all available data provides accurate estimates on the counterfactuals \shi{and} guides a student model which uses fewer features. To avoid overloading notations, in the remainder of the paper, we assume $x^d$ are $d$ categorical features that can be accessed by the student model. Our method can also be applied to numerical features that are discretized, as shown in the experiments in Section \ref{sect_experiments}.  \sun{More discussions can be found in the supplementary material.}

\subsection{Decision Rule Space}\label{sect_rule_space}
The high level idea of our approach is to identify $n$ decision rules from a potentially huge rule space, which contains every possible combination of input features for policy prescription. To model this, we construct  a feature graph, where each decision rule is mapped to a distinct, independent path in the graph. Specifically, we consider an acyclic multi-level digraph, $G(V, E)$, where each feature indicates a level in the graph, represented by multiple nodes corresponding to its distinct feature values. 
Nodes of one feature are fully connected to nodes in the next level. $V$ also contains a  source and sink node. In addition,  a unique and critical characteristic of $G$ is that for each feature, with the exception of the action nodes $\pi$, we introduce a dummy node \emph{SKIP}.  Fig~\ref{feature_space} illustrates an example of a feature graph with three features, including the action nodes.


A decision rule is defined as a path $P\in\mathcal{P}$ from the source to the sink node on $G$. 
\shi{A path passing through a \emph{SKIP} node excludes the corresponding feature from the decision rule.}
As  \emph{SKIP} nodes allow paths to ignore  features, paths on this acyclic graph  represent  all possible feature combinations, defining the full decision rule set $\mathcal{P}$ that one has to enumerate through to solve the policy prescription problem. In the rest of the paper, we use paths and  rules interchangeably. 


\begin{proposition}\label{graph_complexity}
$|V|$ is linear in $d$.  $|\mathcal{P}|$ is exponential in $d$. 
\end{proposition}

The space for decision rules  expands dramatically for high-dimensional datasets and enumerating $\mathcal{P}$ 
%
becomes  intractable. We present an algorithm which searches through the space efficiently in Section~\ref{sect_algorithm}, by solving linear programs (LP) and dynamically generating paths when needed.

So far, we have described our approach of policy prescription as a rule selection problem, whereas the SPMT problem in (\ref{SPMT}) restricts the policy to be a multiway-split tree $\mathcal{T}_m(n)$. There is a key distinction between a set of decision rules and a tree: the former is viewed as unordered rules which may overlap, whereas a tree has a hierarchical structure   which ensures rules are non-overlapping, i.e., each example is covered by exactly one rule (\citealp{Furnkranz2017}). In the next section, we will introduce a MIP formulation to ensure the sample-level non-overlapping condition will be satisfied by the prescribed rules. Note that while the arrangement of features in $G$ is not relevant to the objective value of (\ref{SPMT}), it does affect interpretability. Hence, we utilize insights from the teacher model and arrange the features based on their importance (e.g., \citealp{lundberg2017unified}) when constructing $G$ (see an example presented in the supplementary material).

\subsection{Set Partitioning Problem}
We now present a novel path-based MIP formulation to select a subset of $n$ rules from a feature graph. We denote decision rule as $j = 1,\cdots,N$, where $N=|\mathcal{P}|$. Let $S_j \in [M]$ be the subset of observations which fall into rule $j$. Every sample in $S_j$  is prescribed with a same action $q_j\in \Pi$.  Let $z_j$ be a binary decision variable which indicates whether rule $j$ is selected (1) or not (0), for $j = 1,\cdots,N$. We would like to enforce that each sample $i$ is only assigned to a single rule, and failure to comply with the condition incurs a positive cost $c_i=c$ for all $i$. Define $r_j$ as the expected outcome associated with rule $j$, i.e., $r_j = \sum_{i \in S_j} g_{i,q_j}$. The  \emph{set partitioning problem} ({SPP}) which identifies $n$ out of $N$ decision rules can be written as follows, 

\begin{align}
\quad\max\quad & \sum_{j=1}^N r_jz_j -\sum_{i=1}^M c_is_i\nonumber\\
\textrm{s.t.} \quad & \sum_{j=1}^N a_{ij}z_j + s_i = 1, \quad\forall i=1,\cdots,M\label{constraint_coverage}\\
&\sum_{j=1}^N z_j \leq n \label{capacity}\\
  &z_j\in\{1,0\}, \quad \forall j=1,\cdots,N \nonumber\\&s_i \geq 0, \quad \forall i=1,\cdots,M \nonumber
\end{align}
where $a_{ij}=1$ if sample $i$ satisfies the conditions (features) specified in rule $j$, and 0 otherwise.  
%
While a sample may fall into several rules, \shi{the set partitioning constraint (\ref{constraint_coverage}) along with the non-negative slack variables $s_i$ that are included with a sufficiently large penalty $c_i$, ensure that each sample is ultimately covered by exactly one rule.}
Cardinality constraint (\ref{capacity}) ensures that at most $n$ rules are active in the optimal solution $\mathbf{z}^*$, where $n$ is a user defined input corresponding to $2^k$ leaf nodes in a binary-split tree of depth $k$.

\begin{theorem}\label{them:equivalence}
Active paths in $\mathcal{P}$ which are identified by $\mathbf{z}^*$ correspond to an optimal solution $\mathcal{T}^*_m(n)$ in  SPMT. 
\end{theorem}
We have shown that if we can solve SPP, then we have identified a multiway-split tree with $n$ rules. It is well-known that SPP is NP-Hard (\citealp{wolsey1999integer}). Nevertheless, (near) optimal or high quality feasible solutions can be obtained in practice for moderate sized instances  (\citealp{atamturk1996combined}). Unfortunately, in our case, the cardinality of feasible rules may become prohibitive to solve it directly. 
We will now present an efficient algorithm to overcome this hurdle.

\section*{An Efficient and Scalable Algorithm}\label{sect_algorithm} 
%
The technique of column generation (CG) is used for solving linear problems with a huge number of variables for which it is not practical to explicitly generate  all columns (variables) of the problem matrix. 
More specifically, we consider a \emph{restricted master problem} (RMP) version of SPP, where we  1) consider only a subset of paths, $\hat{N}$, where  $\hat{N}$ is typically much smaller than $N$, and 2) relax the integrality constraints on $z_j$ to $0\leq z_j\leq 1$ for all $j=1,\cdots,\hat{N}$.

Denote the dual variables   associated with the set partitioning constraints in (\ref{constraint_coverage}) and the cardinality constraint in (\ref{capacity})  as ${\lambda}_i$ and $\mu$ respectively. The dual formulation of RMP can be found in the supplementary material. In particular, the dual feasibility constraint corresponding to path $j$ is 
given by

\begin{equation}
    \sum_{i=1}^M a_{ij}\lambda_i+\mu\geq r_j, \quad \forall j=1,\cdots,\hat{N} \label{dual_constraint}
\end{equation}
From (\ref{dual_constraint}), we can derive the reduced cost for path $j$ as
\begin{equation}
    rc_j =  r_j -\left(\sum_{i=1}^Ma_{ij}\lambda_i + \mu\right). \label{eq_reduced cost}
\end{equation}
Additional constraints are handled in (6) similar to the partitioning and cardinality restrictions via 
their corresponding dual variables. 

The key idea of CG is as follows - as we solve a much smaller RMP to optimality, dual feasibility is guaranteed only for the rules included in $\hat{N}$ and we must verify that the dual feasibility are also satisfied by the rules not included in the RMP. A path violating the dual constraint  (\ref{dual_constraint}) has a positive reduced cost and must be added to the RMP. 
Therefore, we identify paths that maximally violate Eq (\ref{dual_constraint}), i.e., $\max rc_j $ or $\min -rc_j$. 
The implicit search for $K$ paths having the highest reduced cost amounts to optimizing a \emph{subproblem}, which is a 
\emph{K-shortest path problem} (KSP) over $G$  (\citealp{horne1980finding, irnich2005shortest}), where the cost on a path is the negative of its reduced cost defined in (\ref{eq_reduced cost}). 
The CG procedure converges to an optimal solution of the LP relaxation of SPP when dual feasibility is achieved, i.e., $rc_j \leq 0$  for all $j = 1,\cdots,N$. Otherwise, new paths are found and their corresponding columns $A_j$ which comprise of coefficients $a_{ij}$  are added to RMP, which is re-optimized. We repeatedly solve RMP and KSP until successive dual solutions converge  or we reach a maximum iteration limit.

Once the CG procedure has converged, we reimpose the binary restrictions on $\mathbf{z}$ and solve the resultant Master-MIP. Any of two MIP approaches below can be adopted \shi{either individually or in combination} depending on the end goal. Since SPP is NP-Hard, obtaining a provably optimal solution may require inordinate run times for challenging instances. For simplicity and replicability, we solve the Master-MIP directly here using a standard optimization library (\citealp{cplex2020}) to distill a (near) optimal subset of prescriptive decision rules. 
\shi{Thereafter, one can optionally switch to} the more advanced ``branch-and-price'' method 
%
\shi{to converge to an} optimal MIP solution \cite{barnhart1998branch}.

\sun{One common yet critical requirement of prescriptive analytics is to incorporate a myriad of operational constraints. They}
can be broadly categorized into  \emph{intra-rule} and \emph{inter-rule constraints}, which are handled differently in the algorithm. 

The scope of \emph{intra-rule constraints} is limited to a single decision rule. They may involve complex path-dependent nonlinear conditions involving several features that cannot be abstracted efficiently into linear constraints.
For instance, disallowing certain feature combinations or  limiting the number of features. 
These constraints are easily handled within the KSP subproblem as a feasibility check while extending a partial path to the next node in $G$ (we refer readers to Section 1.2 in \citealp{barnhart1998branch} for more discussions). 

\emph{Inter-rule constraints} span across many decision rules. One example is the cardinality constraint in (\ref{capacity}). 
Inter-rule constraints can be expressed as linear inequalities  in the RMP, which in turn influence the subproblem \shi{by} modifying the reduced cost. 
We present detailed use cases in Section~\ref{sect_real} along with a house pricing constraint and loyalty pricing implementation in the supplementary material.

\begin{figure*}
\begin{subfigure}{.5\textwidth}
  \includegraphics[width=.85\linewidth]{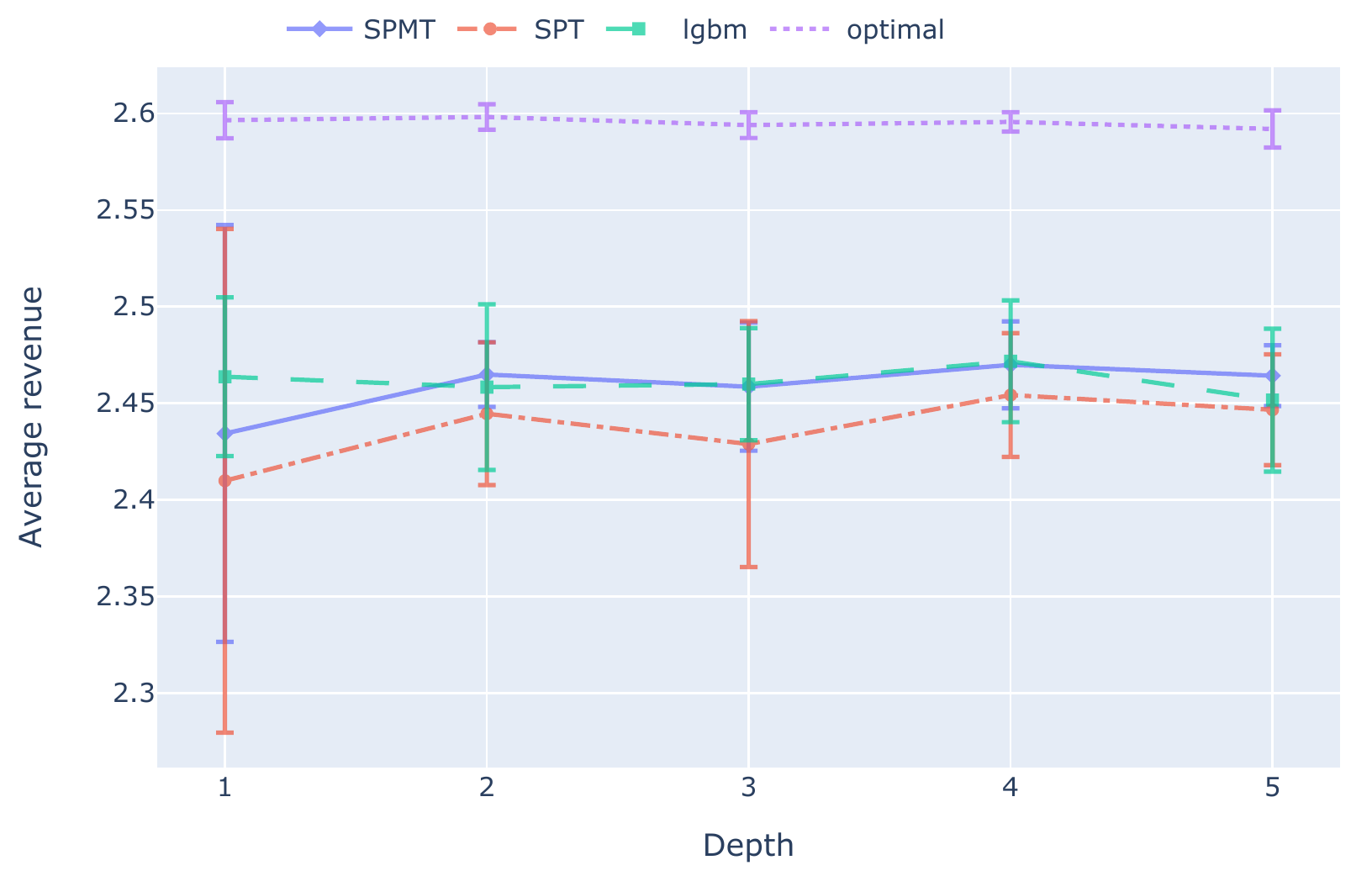}
  \caption{Revenue with varying tree depth}
  \label{rev_with_depth}
  \end{subfigure}
\begin{subfigure}{.5\textwidth}
  \includegraphics[width=.85\linewidth]{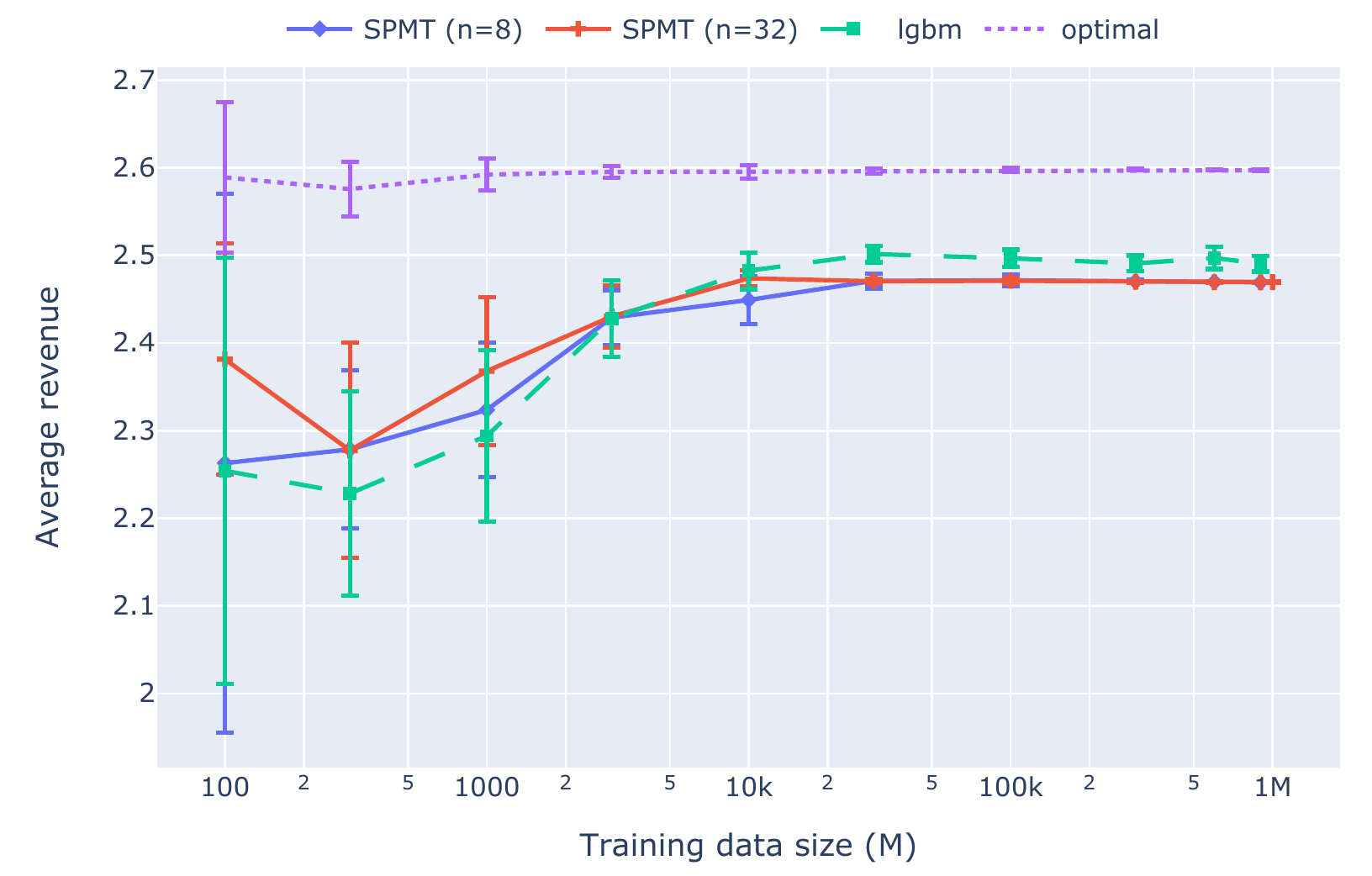}
  \caption{Revenue with varying training data size}
  \label{rev_with_datasize}
  \end{subfigure}
  \caption{Experiment results on realized outcome for synthetic Dataset (6)}
\end{figure*}

\section*{Experiments}\label{sect_experiments}

In this section, we perform experiments to demonstrate the efficacy of our method, \emph{Student Prescriptive Multiway-split Tree} (SPMT).  We first perform experiments on six simulated datasets which allow us to accurately calculate the resulting outcome. We want to point out that since none of the existing interpretable tree-based methods are capable \shi{of} handling constraints, we can only benchmark against them in the \emph{unconstrained} setting.
Next, we present three use cases with real data as we discuss constrained policy prescription, two of which use publicly available data. 

The following baselines are considered: The greedy student prescriptive tree (\textsc{SPT}), is used as the key tree-based prescriptive benchmark since \citealp{biggs2020model} has demonstrated that it outperforms other competing prescriptive tree  methods (i.e., \citealp{kallus2017recursive,athey2016recursive}). 
A true optimum policy (\textsc{optimal})  which can be found by identifying the  price which results in the highest revenue according to the underlying probability model in synthetic datasets. The baseline (\textsc{lgbm}) determines the sample-level optimal policy based on the predictions of a  lightGBM  model, which is also the teacher model. 
\shi{All experiments were run on an Intel 8-core i7 PC with 32GB RAM.} 
CPLEX 20.1  was used to solve the RMP.  
The minimum number of samples per rule for SPMT is set to 10. $K=100$ for KSP. 

\subsection{Synthetic Datasets}\label{sect_synthetic_data}
We follow the identical setup in \citealp{biggs2020model}, where we generate six datasets consisting entirely of numerical features, simulated from different generative demand models. 
The data generation process can be found in the supplementary materials. We train a gradient boosted tree ensemble model using the lightGBM package (\citealp{ke2017lightgbm}) as the teacher model. 
The action set consists of 9 prices, ranging from the $10^{th}$ to the $90^{th}$ percentile of observed prices in $10\%$ increments.


\noindent\textbf{Varying policy size} We explore how the expected revenue changes with the number of policies. We use a SPT of depth $k$ which is grown to its full width to  benchmark against a  multiway-split tree (SPMT) with $n=2^k$ rules, where $k=\{1,..,5\}$. 
The number of training samples is held constant at $M=5000$.  
For each tree depth, we run 10 independent simulations for each dataset and each $k$ (i.e., 300 experiments in total). We highlight the results for $k=3$ and 5 in Table~\ref{revenue_table}, which summarizes mean revenue over the simulations. 
Fig~\ref{rev_with_depth} shows the performance on dataset (6) with respect to the number of rules as an illustrative example, while detailed results for individual datasets are included in the supplementary materials.

%
\shi{Although both  SPMT  and  SPT  use  the same  teacher model, SPMT generally outperforms SPT, which employs a greedy heuristic, as seen in the modest gains shown in the experiments. 
%
Furthermore,} SPT is unable to incorporate constraints, a shortcoming which severely  limits  its  usefulness.  

\begin{table*}[t]
\makebox[\textwidth][c]{
    \begin{tabular}{cccccccc}
\toprule
      & &      & \multicolumn{2}{c}{n = 8 / k = 3} & \multicolumn{2}{c}{n = 32 / k = 5 } \\
\cmidrule(lr){4-5}\cmidrule(lr){6-7}
 Dataset & optimal  & lgbm & SPMT & SPT   & SPMT & SPT \\ \midrule
    1 & 3.275 (0.016) &3.161 (0.054) & {3.233} (0.038) & 3.201 (0.071) & {3.201} (0.053) & 3.188 (0.063) \\
      2 & 2.963 (0.034) & 2.716 (0.038) & {2.166} (0.034) & 2.127 (0.042) & {2.378} (0.015) & 2.277 (0.039) \\
      3 &  3.420(0.006) & 3.275 (0.064) & {3.347}(0.048) & 3.331(0.052) & {3.336}(0.047) & 3.333(0.047) \\
      4 & 3.492(0.012) & 3.344 (0.057) & {3.365}(0.051) & 3.343(0.079) & {3.379}(0.055) &  3.366(0.06) \\
      5 &  3.35(0.007) & 3.260 (0.029) & {3.291}(0.056) &  3.260(0.067) & {3.296}(0.017) & 3.287(0.022) \\
      6 &  2.592(0.01) & 2.456 (0.033) & {2.459}(0.033) & 2.429(0.064) & {2.464}(0.016) & 2.447(0.029) \\
\bottomrule
\end{tabular}%
}
\caption{Comparison on realized outcome for synthetic datasets ($M=5000$)}
\label{revenue_table}\end{table*}%

\noindent\textbf{Varying data size} We also investigate how SPMT performs as the size of the samples increases, where  
$M=\{100, 300, 10^3, 3\times10^3, 10^4, 3\times10^4, 10^5, 3\times10^5, 6\times 10^5, 9\times 10^5, 10^6\}$.
We focus on $n=8$ and 32 rules, again for 10 independent iterations for each  configuration (i.e., 220 experiments) using Dataset (6). Fig~\ref{rev_with_datasize} shows the result on realized revenue versus data size $M$. \sun{\footnote{We omitted SPT  here because its current implementation from \citealp{biggs2020model} is struggling with datasets larger than 30K samples. See supplementary material for more discussions. }}

The baseline \textsc{lgbm}  is an upper bound to  SPMT  in theory since the latter trades off performance for interpretability. It can be seen as the revenue gap between \textsc{lgbm} and SPMT policies Fig~\ref{rev_with_datasize} when the training sample size becomes large and the lgbm model produces a more accurate estimate of the underlying groundtruth. With small sample sizes, \textsc{lgbm} may overfit and the realized outcome may appear worse/noisier. It is crucial to note that policy from the \textsc{lgbm} baseline is not interpretable and can be potentially very complex (such as a fully personalized price). 
It is also not surprising that a policy with more rules generally outperforms fewer rules. However, once the sample size exceeds $3\times 10^4$, both policies with $n=8$ and $32$  converge to the same realized revenue and we do not observe further improvement from either policy by increasing the sample size. 

\noindent\textbf{Runtime study} We investigate CG runtime  as the sample size $M$ increases from 100 to $10^6$, with the same setup described earlier. 
In our experiments, all \shi{RMP} instances \shi{converged in less than 20 CG iterations}, with 11.19 (1.33) seconds for each iteration of RMP and KSP. Detailed results are available in the supplementary material.
The key observation is that our \shi{proposed SPMT approach} is inherently more scalable than prior MIP approaches \shi {when near-optimal solutions suffice}.
In the experiment which allowed \shi{up to 32 rules for the dataset having $10^6$ samples}, prior MIP models would require more than $32\times10^6$ binary variables. In comparison, our Master-MIP required less than 3000 binary variables.

\subsection{Real Datasets}\label{sect_real}
\noindent\textbf{Grocery pricing case study}
We consider the problem of grocery pricing, using a publicly available dataset\footnote{\url{https://www.dunnhumby.com/careers/engineering/sourcefiles}} that contains over two years of household level transactions from a group of 2,500 households  at a grocery retailer.  A processed dataset with 97,295 rows which contain both purchases and ``no-purchases'' of strawberries is available online.\footnote{\url{https://docs.interpretable.ai/stable/examples/grocery_pricing/}}

The unconstrained use case of maximizing expected revenue has been studied in \citealp{biggs2020model} and \citealp{amram2020optimal}. 
We follow the same experiment setup described in these papers, where we divide the data in halves, and independently train a teacher and an evaluator (as the ground truth) using lightGBM. With $n=64$, SPMT policy achieves 82.8\% increase in revenue over the historical baseline, compared to 65.9\% and 77.1\% reported in \citealp{biggs2020model} and \citealp{amram2020optimal} respectively. Noting that the latter work and our approach both \shi{seek optimal prescriptive policies}, the difference in gain highlights the benefit of a multiway-split tree over a binary-split tree, for being more informative. An example of a SPMT with $n=8$ is shown in Fig~\ref{multiway_split}. 

There are two drawbacks with the existing unconstrained policy. First, its predicted demand almost doubles the historical sales, which may not be realizable due to supply constraints and other limitations. Second, the price discrimination based on customer-level features raises legal and fairness concerns.  
We briefly describe two constrained scenarios to address these issues, where more details are provided in the supplementary material. First, stores are grouped into 5 clusters based on historical sales and we require the predicted demand under a new pricing policy does not deviate too much from the historical average. With $n=64$, SPMT produces rules that result in 47.49\% and 61.72\% increase in revenue when the allowable change in predicted demand is set to 25\% and 50\% of the historical value. In a separate scenario, we use inter-rule constraints to model a store-based loyalty pricing policy, where all shoppers at a store receive the same price except loyalty-card shoppers, who receive a price that is same or lower. \shi{This setting is an example of inter-rule constraints involving logical conjunctions of categorical and numerical features (StoreID, Loyalty, Price), which are handled naturally in our path-based MIP, but relatively difficult to express efficiently using arc-based feature-level models.} The resulting \shi{pricing} policy using only storeID and a loyalty indicator to define rules achieves 65.4\% increase in revenue  over the baseline without discriminating customers based on their personal features. 

%


\noindent\textbf{{Housing upgrade use case}}
We consider a Kaggle dataset on house prices\footnote{\url{https://www.kaggle.com/harlfoxem/housesalesprediction}} that includes numerical and categorical attributes such as age, square footage, number of bedrooms, bathrooms, zipcode, as well as grade, i.e., an index from 1-13, where higher value indicates better construction and design.    
We first train a regression model to predict the house value given the  attributes as the teacher. We then formulate a prescriptive problem, i.e., maximize the predicted sale price in the test set by adjusting the grade of individual houses. As the teacher model predicts that grade increases the house value, an unconstrained policy search simply recommends that all houses upgrade to the highest level. 

One of the constrained scenarios that we consider requires houses to ``conform to the neighborhood'', i.e., the predicted sale price of houses with the same zip code must be within 10\% of the historical value. The decision rules no longer upgrade all houses but prescribe a mix of upgrades and downgrades based on other housing features while satisfying the price constraint for each neighborhood. More constrained scenarios are included in the supplementary materials.

\noindent\textbf{Airline pricing use case}
We consider the task of pricing premium seats for a large airline \sun{(an industry partner)} using easily interpretable prescriptive rules. 
%
%
The dataset contains more than 1.3 million samples across multiple markets with several  features such as departure day and time, duration of the flight, market, whether one-way or round trip. \sun{We want to point out that our work allows  multiple actions to be added as stacks of action nodes
%
e.g., to represent the prices for the first class, business class and premium economy seats respectively. Structured decision making can be modeled as constraints, e.g., given the same trip information, a business class (premium economy) seat should be at least \$250 cheaper than a first (business) class seat.} For confidentiality reasons, details of the study are omitted.

We compare the policies from the greedy SPT and our proposed SPMT method for a same number of decision rules. 
With SPT, seat capacity limits were violated and resulted in a overly high predicted gain. Further analysis revealed that these rules would have allowed early bookers to purchase premium seats at a bargain,  cannibalizing the demand from high-value customers who typically show up later. By adding market-specific capacity constraints on the predicted conversions, we revised the rules to achieve a more realistic gain that is more likely to be realized in live settings.

\section*{Conclusion}
Our work adds to a stream of  growing research  which deserves more attention in the literature due to its paramount practical importance, i.e., generating interpretable policy from data for decision making.
We address the problem of distilling prescriptive rules using a teacher model, providing counterfactual estimates. We formulate a mixed integer optimization problem to distill active decision rules that are shown to represent a multiway-split tree on an acyclic feature graph. 
The number of binary decision variables in our model 
can increase exponentially in the number of features that comprise the rules. As only a small subset of rules are active in an optimal solution, we propose a 
column generation approach to solve the MIP. Our method can handle complex inter-rule and intra-rule constraints that cannot be satisfied by existing prescriptive tree methods. We provide computational results on synthetic and publicly available real-life data and solve  instances having up to a million samples within an hour of computation time.
The proposed method offers a generic and widely applicable way of generating (near) optimal prescriptive rules by transforming blackbox AI/ML predictions into operationally effective decisions that may benefit enterprises. 
    



\bibliography{ref}

\bibliographystyle{plainnat}
\bibliography{ref}

\end{document}


\maketitle

\normalsize
\section{Section 3 - Problem Formulation}
\subsection{Proof of Proposition 1} 
In the acyclic multi-level directed graph which we consider, $G(V, E)$,  each feature indicates a level in the graph, represented by multiple nodes corresponding to its distinct feature values. 
For each feature, with the exception of the action nodes $\pi$, we introduce a dummy node \emph{SKIP}. Without the loss of generality, with a $d$-dimensional input vector $x$, action feature $\pi$ is the last feature (see Fig 1a as an illustrative example). Denote $L_f$ as the number of unique feature values for the $f${-th} feature. 

The total number of nodes in $G$ is given by $|V|=\sum_{f=1}^{d-1}(L_f+1) + L_d + 2$, including the source and sink nodes, i.e., $|V|=O\left(\sum_{f=1}^d L_f\right)=O\left(d\right)$. On the other hand, the number of paths from the source to the sink is given by $|\mathcal{P}|=O\left(\prod_{f=1}^d L_f\right)=O\left(L^d\right)$, for some constant $L$.

\subsection{Proof of Theorem 1}
By including \emph{SKIP} nodes, the path set $\mathcal{P}$ which corresponds to the decision rule set includes every possible combination of features. Moreover, the set coverage constraint (3) in SPP has ensured that rules are non-overlapping, i.e., each sample can only be assigned to a single rule. Then we reach the desired result.

\subsection{Remarks on numerical features} 
For tree-based approaches, numerical input are typically discretized. For example, consider a numerical feature with values in [0, 1] discretized into $\kappa=3$ levels, corresponding to 3 intervals [0, 0.33), [0.33, 0.67) and [0.67, 1.0]. A naive approach is to create 3 nodes in a multiway-split tree. However, this approach is limiting, considering a binary-split tree can branch on conditions such as $x\leq 0.67$ or $x> 0.3$. We consider a \emph{cumulative binning} method, where intervals can be overlapping. More specifically, we create additional nodes {[0, 0.67), [0, 1.0], [0.33, 1.0]}, yielding a total of 6 nodes. 
 Note that the MIP formulation includes a coverage constraint to ensure that the final rules do not contain overlapping samples. More generally, this binning method results in \textit{O}($\kappa^2$) nodes, i.e, an improvement in the expressiveness of a rule at the expense of higher computational complexity. More efficient methods to handle numerical features can be employed as a topic for future work.

\shi{A similar cumulative binning approach can be adopted to manage categorical features if we want to create nodes associated with multiple levels of a feature. Doing so can improve solution quality at the expense of increased rule complexity.}

\section{Section 4 - An Efficient and Scalable Algorithm}

\textbf{Dual formulation of RMP}

\noindent For the \emph{restricted master problem} (RMP),  we  1) consider only a subset of paths, $\hat{N}$, where  $\hat{N}$ is typically much smaller than $N$, and 2) relax the integrality constraints on $z_j$ to $0\leq z_j\leq 1$ for all $j=1,\cdots,\hat{N}$. 
Denote the dual variables   associated with the set partitioning constraints in (3) and the cardinality constraint in (4)  as ${\lambda}_i$ and $\mu$ respectively. The dual of the RMP can be written as the following,  
\begin{align*}
\textbf{(Dual of RMP)} \quad\max\quad & \sum_{i=1}^M \lambda_i + n\mu\nonumber\\
\textrm{s.t.} \quad & \sum_{i=1}^M a_{ij}\lambda_i+\mu\geq r_j, \quad \forall j=1,\cdots,\hat{N}\label{dual_constraint}\\
& \lambda_i\geq -c_i, \forall i = 1,\cdots, M\nonumber\\
& \mu\geq0 \nonumber
\end{align*}
Note that since $z_j\leq 1$ is implied by Eq (3), we are only left with $z_j \geq 0$ in the primal RMP problem. 

\textbf{Remark:} \sun{A feasible solution to RMP is to set all the slack variables $s_i$ to 0}
%
%
and add a column that prescribes a single rule for all data samples.
The corresponding path in the feature graph passes (only) through all SKIP nodes and then through any of the action nodes.

\section{Section 5 - Experiments}
\subsection{Synthetic datasets}
\subsubsection{Data generation}
We follow the identical setup in \cite{biggs2020model}, where six datasets consisting entirely of numerical features are  simulated from different generative demand models. For completeness, we include the data generation process below. 

The datasets we examine come from 
the following generative model:
\begin{equation}
	Y^*=g(X)+h(X)P + \epsilon, \quad Y=
    \begin{cases}
      1, & \text{if}\ Y^*>0 \\
      0, & \text{if}\ Y^*\leq 0
    \end{cases} \label{latent_eq} \\ 
\end{equation}
\begin{itemize}
\item Dataset 1: linear probit model with no confounding: $g(X)=X_0$, $h(X)=-1$ and $X \sim N(5,I_2)$ and $P \sim N(5,1)$.
\item Dataset 2: higher dimension probit model with sparse linear interaction $g(X)=5,~ h(X)= - 1.5 (X'\beta)$,  $\{X_i\}_{i=1}^{20} , \{\beta_i\}_{i=1}^5, \epsilon_i \sim N(0,1), P_i\sim N(0,2), \{\beta_i\}_{i=6}^{20} = 0$, where the purchase probability is only dependent on the first 5 features.
\item Dataset 3: probit model with step interaction: $g(X)=5$,  $h(X)=-1.2 \mathbbm{1}\{X_0< -1\} -1.1 \mathbbm{1}\{-1 \leq X_0< 0\}  -0.9 \mathbbm{1}\{0 \leq X_0< 1\} -0.8 \mathbbm{1}\{1 \leq X_0 \}$.
\item Dataset 4: probit model with multi-dimensional step interaction: $g(X)=5$, $h(X)=-1.25 \mathbbm{1}\{X_0< -1\} -1.1 \mathbbm{1}\{-1 \leq X_0< 0\}  -0.9 \mathbbm{1}\{0 \leq X_0< 1\} -0.75 \mathbbm{1}\{1 \leq X_0 \} - 0.1\mathbbm{1}\{X_1< 0\}+  0.1\mathbbm{1}\{X_1 \geq 0\} $.

\item Dataset 5: linear probit model with confounding: $g(X)=X_0$, $h(X)=-1$, $X \sim N(5,I_2)$.
\item Dataset 6: probit model with non-linear interaction: $g(X)= 4|X_0+X_1|$,     $h(X)=-|X_0+X_1|$.
\end{itemize}

\begin{figure}
  \begin{subfigure}{.5\textwidth}
  \includegraphics[width=1.1\linewidth]{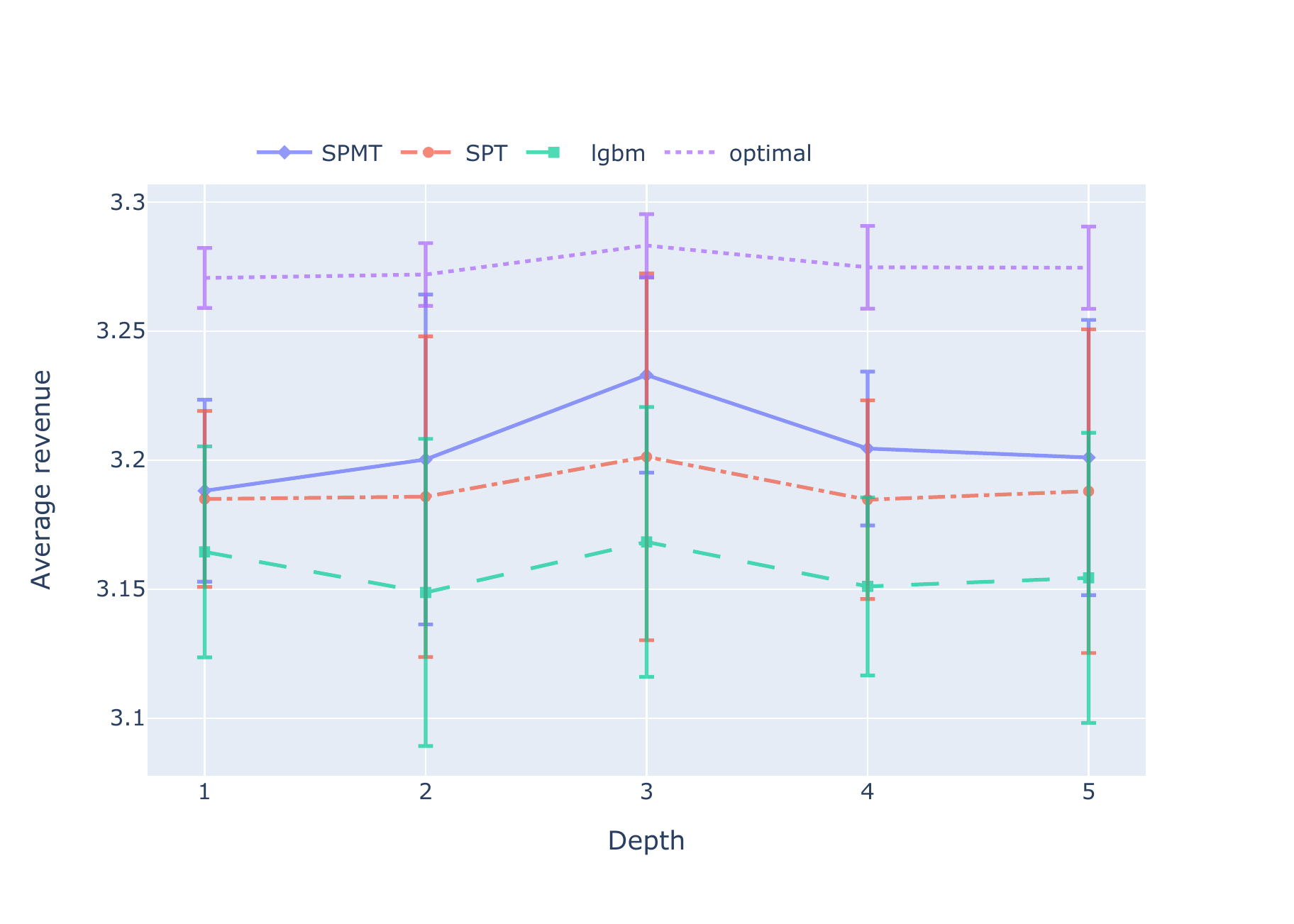}
  \caption{Dataset 1}
  \end{subfigure}
\begin{subfigure}{.5\textwidth}
  \includegraphics[width=1.1\linewidth]{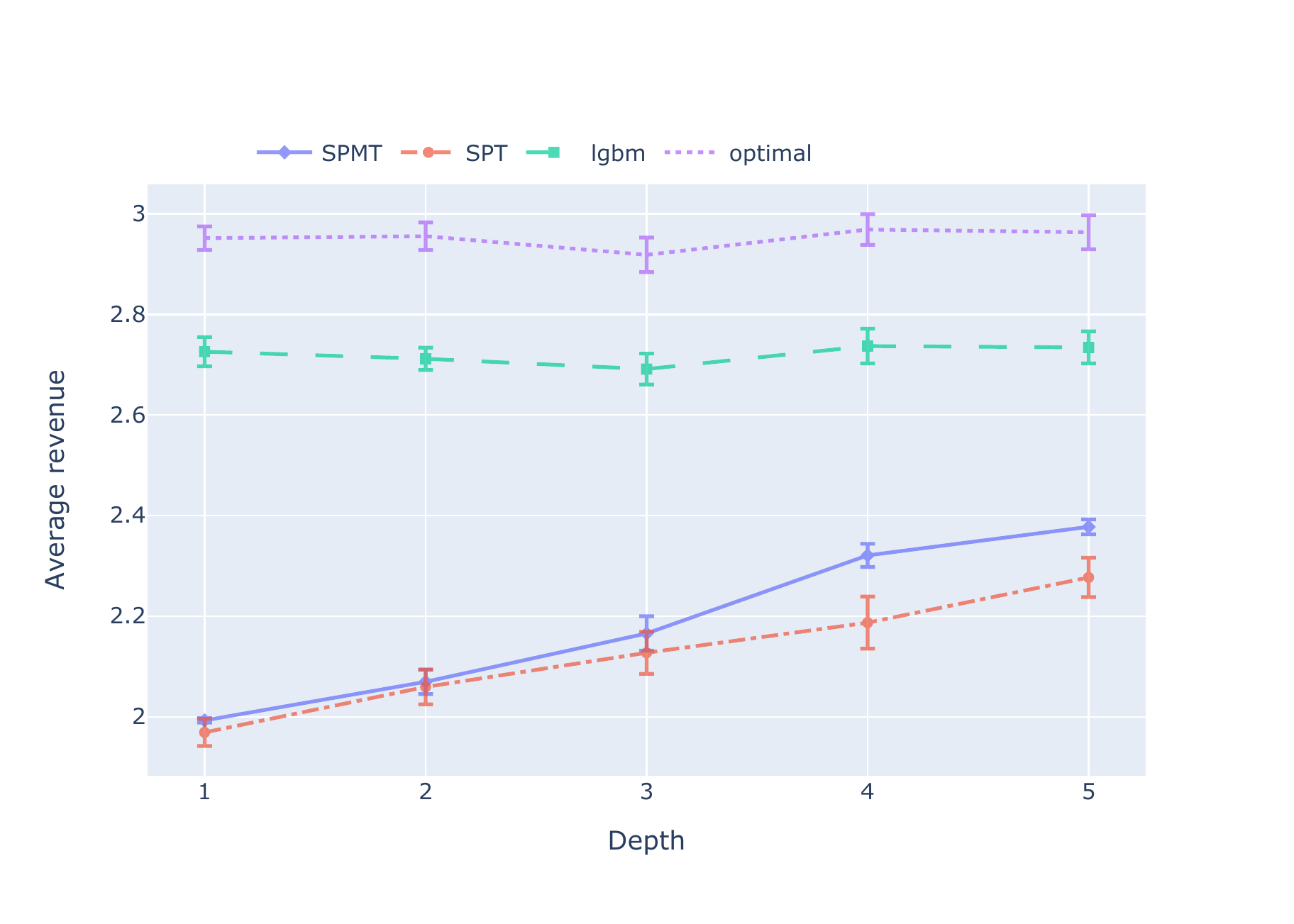}
  \caption{Dataset 2}
  \end{subfigure}
\begin{subfigure}{.5\textwidth}
  \includegraphics[width=1.1\linewidth]{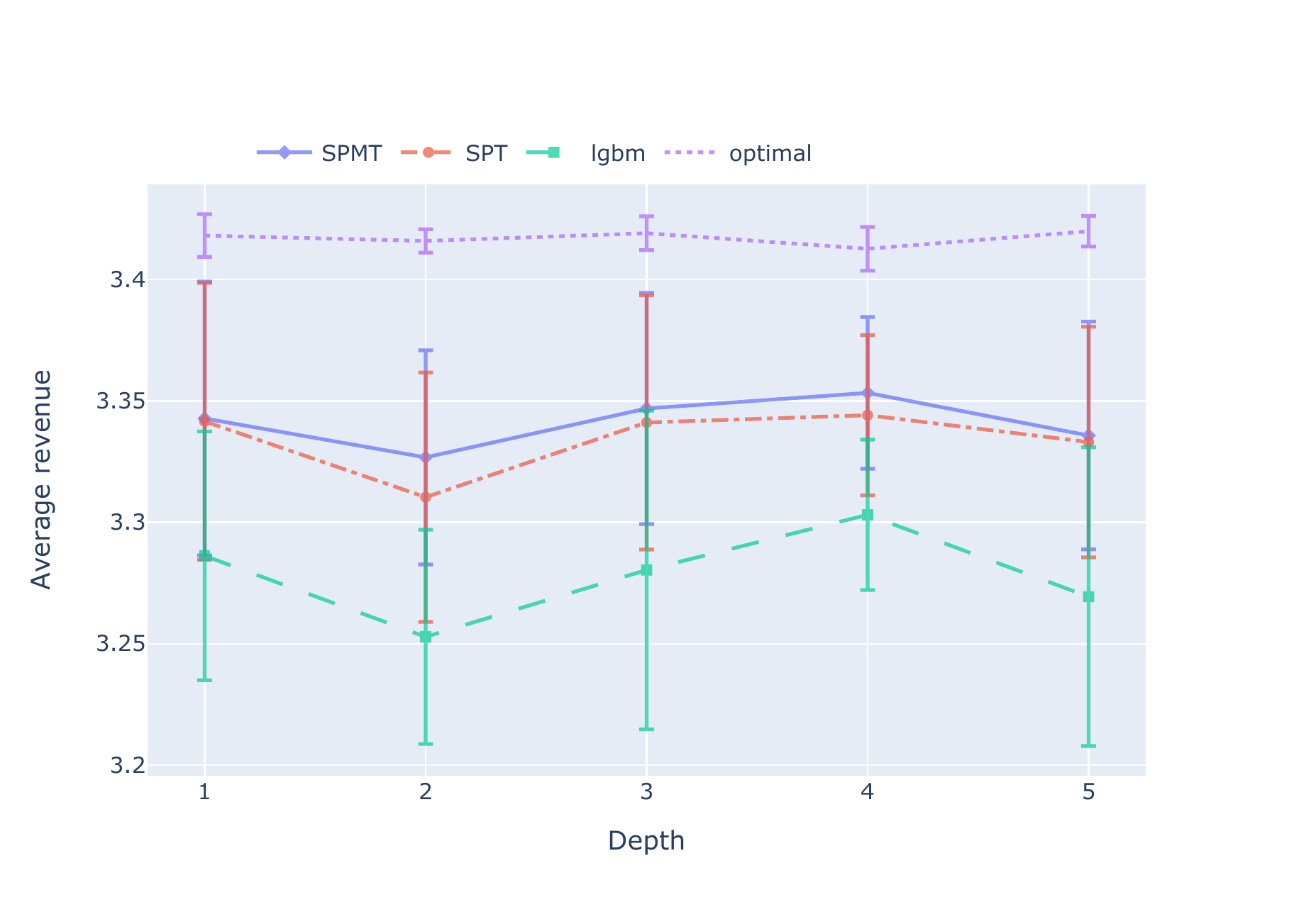}
  \caption{Dataset 3}
  \end{subfigure}
  \begin{subfigure}{.5\textwidth}
  \includegraphics[width=1.1\linewidth]{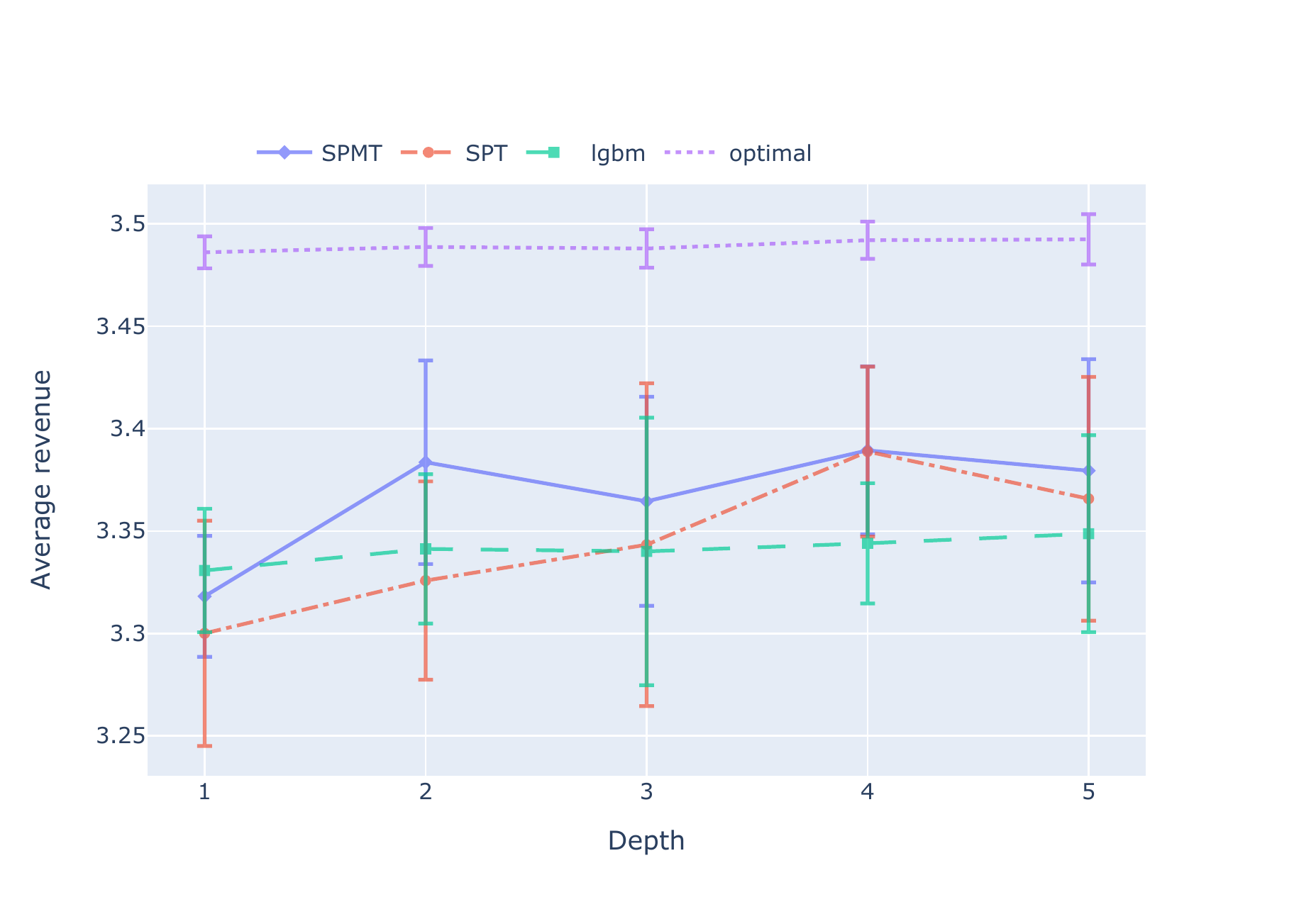}
  \caption{Dataset 4}
  \label{depth_datset_4}
  \end{subfigure}
 \begin{subfigure}{.5\textwidth}
  \includegraphics[width=1.1\linewidth]{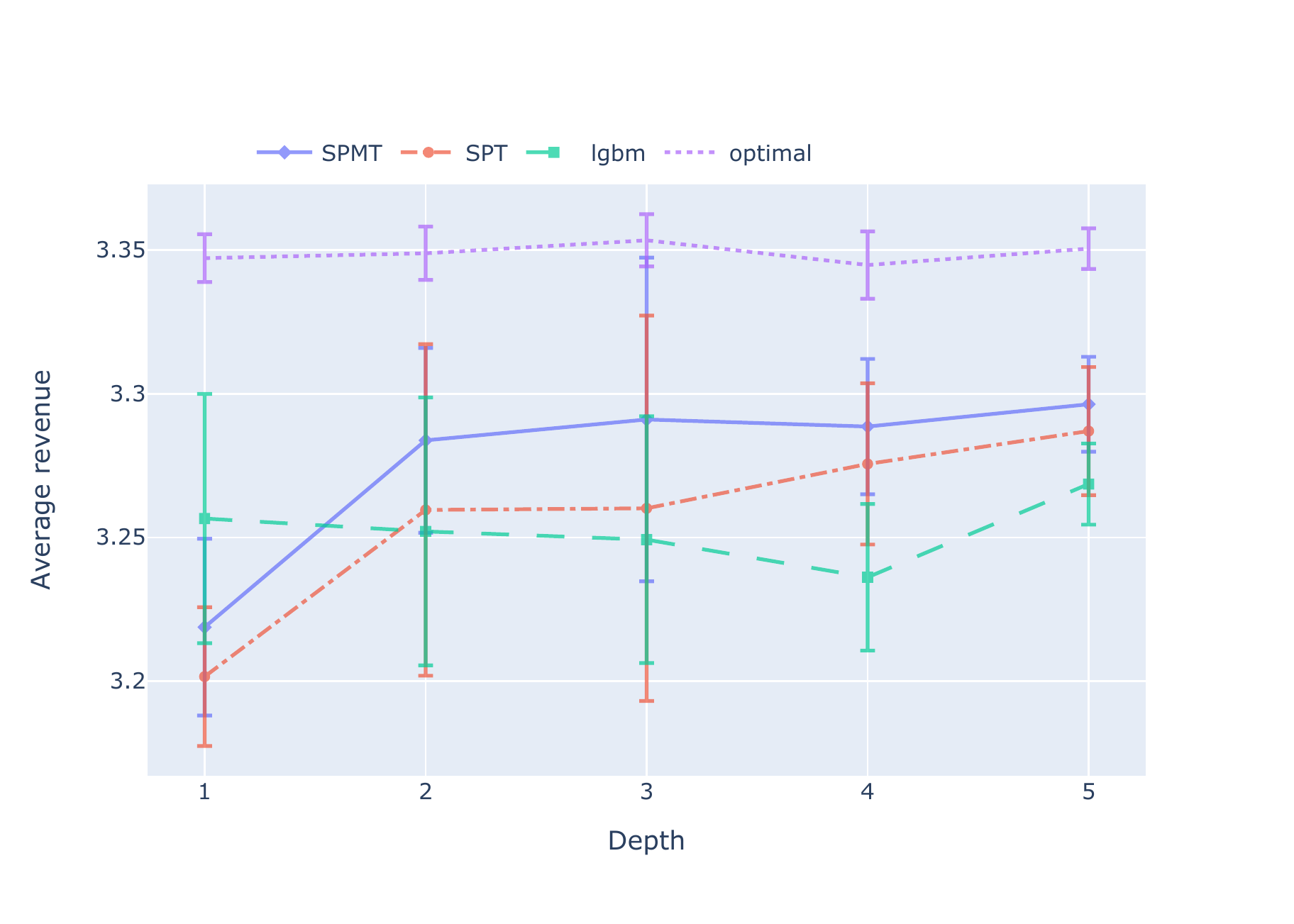}
  \caption{Dataset 5}
  \end{subfigure}
  \begin{subfigure}{.5\textwidth}
  \includegraphics[width=1.1\linewidth]{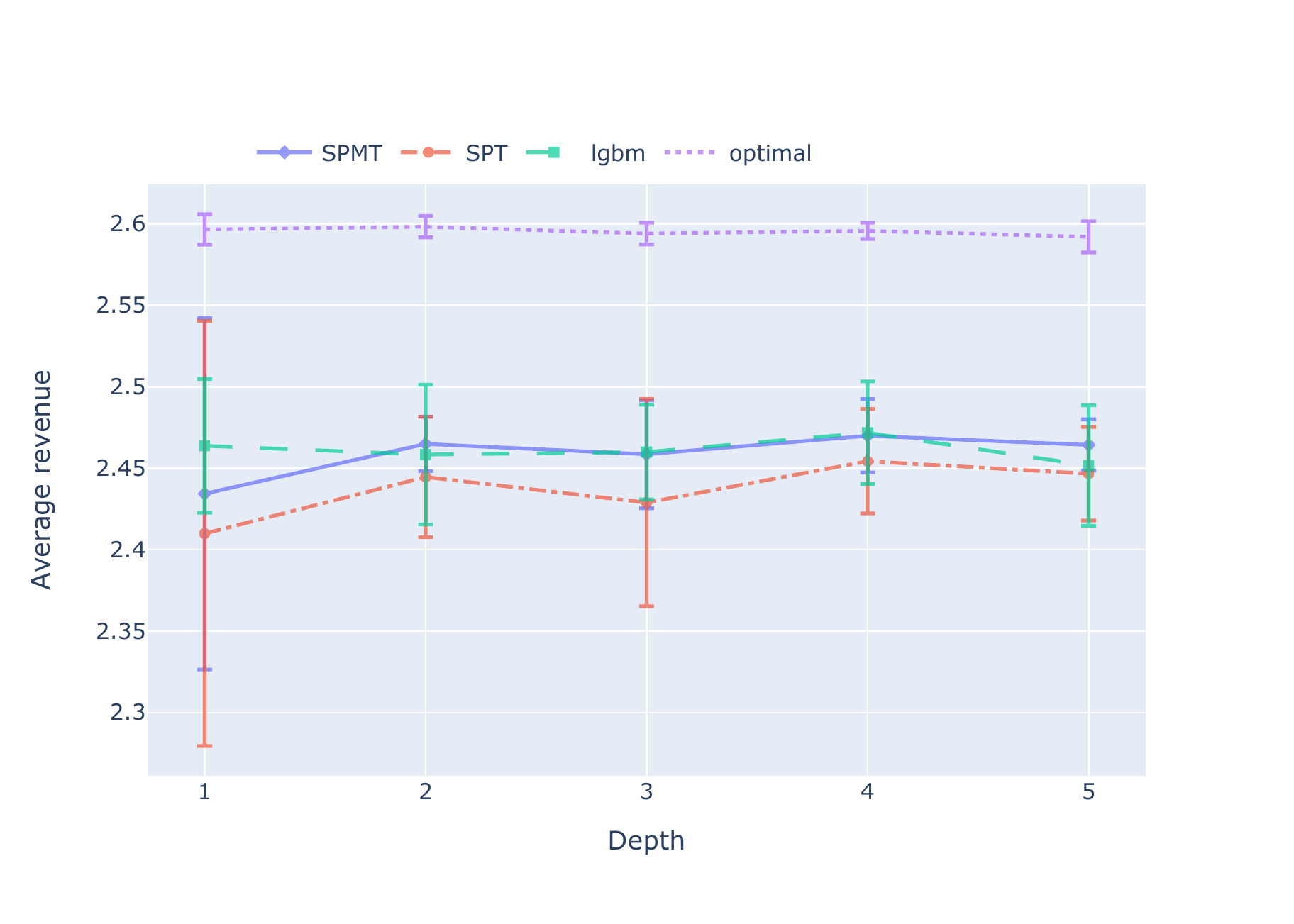}
  \caption{Dataset 6}
  \end{subfigure}
  \caption{Experiments synthetic data with increasing tree depth $k$ of SPT, equivalently rule complexity $n=2^k$ in SPMT}
  \label{synth_vary_depth}
\end{figure}

We set $X= (X_0, X_1) \sim N(0,I_2)$, $P \sim N(X_0+5,2)$,   $\epsilon \sim N(0,1)$, i.i.d. unless otherwise mentioned.  Dataset (1) and (5) are a common demand model used in the pricing literature, with and without confounding effects.  Datasets (2)-(4) and (6), have heterogeneity in treatment effects making them suitable for testing personalized pricing.  All datasets except dataset (1) and (2) have confounding of observed features, where the price observed is dependent on the features of the item.
For the teacher model, we train a gradient boosted tree ensemble model using the lightGBM package \cite{ke2017lightgbm}. We use default parameter values, with 50 boosting rounds. 
The discretized price set is set to 9 prices, ranging from the $10^{th}$ to the $90^{th}$ percentile of observed prices in $10\%$ increments. 

\subsection{Additional experiment results}

\subsubsection{Varying policy size and data size}
Figure \ref{synth_vary_depth} shows for each of the six datasets, how the expected revenue changes with respect to  the depth of the SPT $k$, or equivalently the rule cardinality in SPMT where $n=2^k$, for the baselines we benchmark against.

Note that SPT is not shown in Fig 2b) in the main paper, as the code used to produce SPT (courtesy of the authors, \cite{biggs2020model}) was struggling to handle datasets larger than 10K samples and taking several hours to run for each iteration. As a result, we skip this baseline for the experiment in Fig 3b) where we vary sample sizes to $10^6$. 

\subsubsection{Run time study}
We investigate column generation (CG) runtime as the sample size $M$ increases, where 
$M=\{100, 300, 10^3, 3\times10^3, 10^4, 3\times10^4, 10^5, 3\times10^5, 6\times 10^5, 9\times 10^5, 10^6\}$.
We focus on $n=8$ and 32 rules, again for 10 independent iterations for each  configuration (i.e., 220 experiments) using Dataset (6).

The run time per iteration of CG equals the sum of the RMP and KSP run times. The RMP, being a linear program can be solved in polynomial time (\cite{bazaraa2008linear}). The KSP is solved using a naive implementation of a $K$-shortest path search on an acyclic graph and is bounded by \textit{O}($K|V|^2M~log~M$). The first term ($K|V|^2$) represents the worst case run time for a $K$-shortest path algorithm on an acyclic graph (\cite{horne1980finding}), and the $(M\log M)$ term captures the time required to update the training samples (sort and store as trees) covered by the path after every partial path extension. More efficient implementations can employed to reduce the KSP run time.

\begin{figure}
\begin{subfigure}{.5\textwidth}
  \includegraphics[width=1\linewidth]{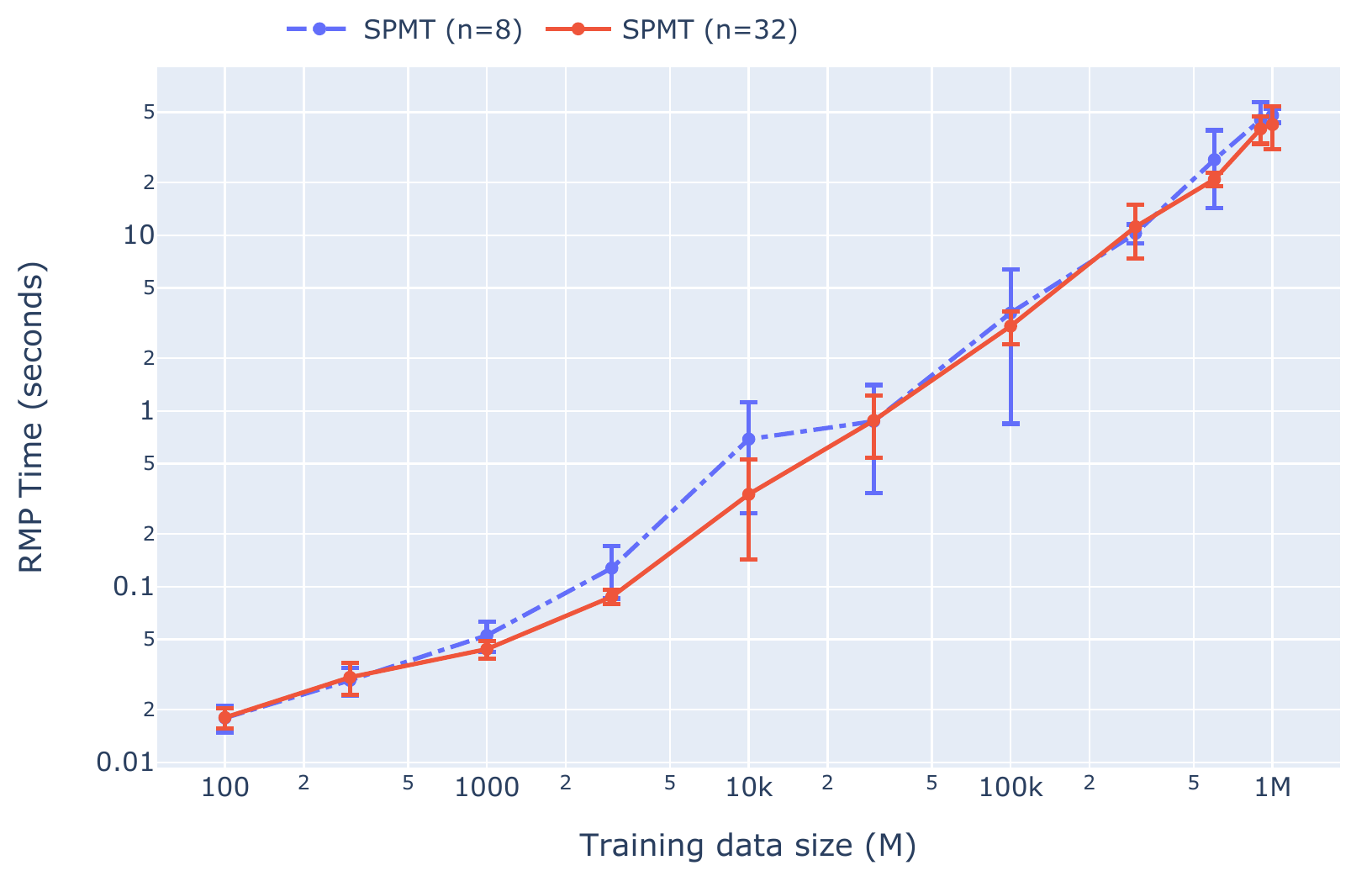}
  \caption{Average  RMP duration}
  \label{rmp_time}
  \end{subfigure}
\begin{subfigure}{.5\textwidth}
  \includegraphics[width=1\linewidth, clip]{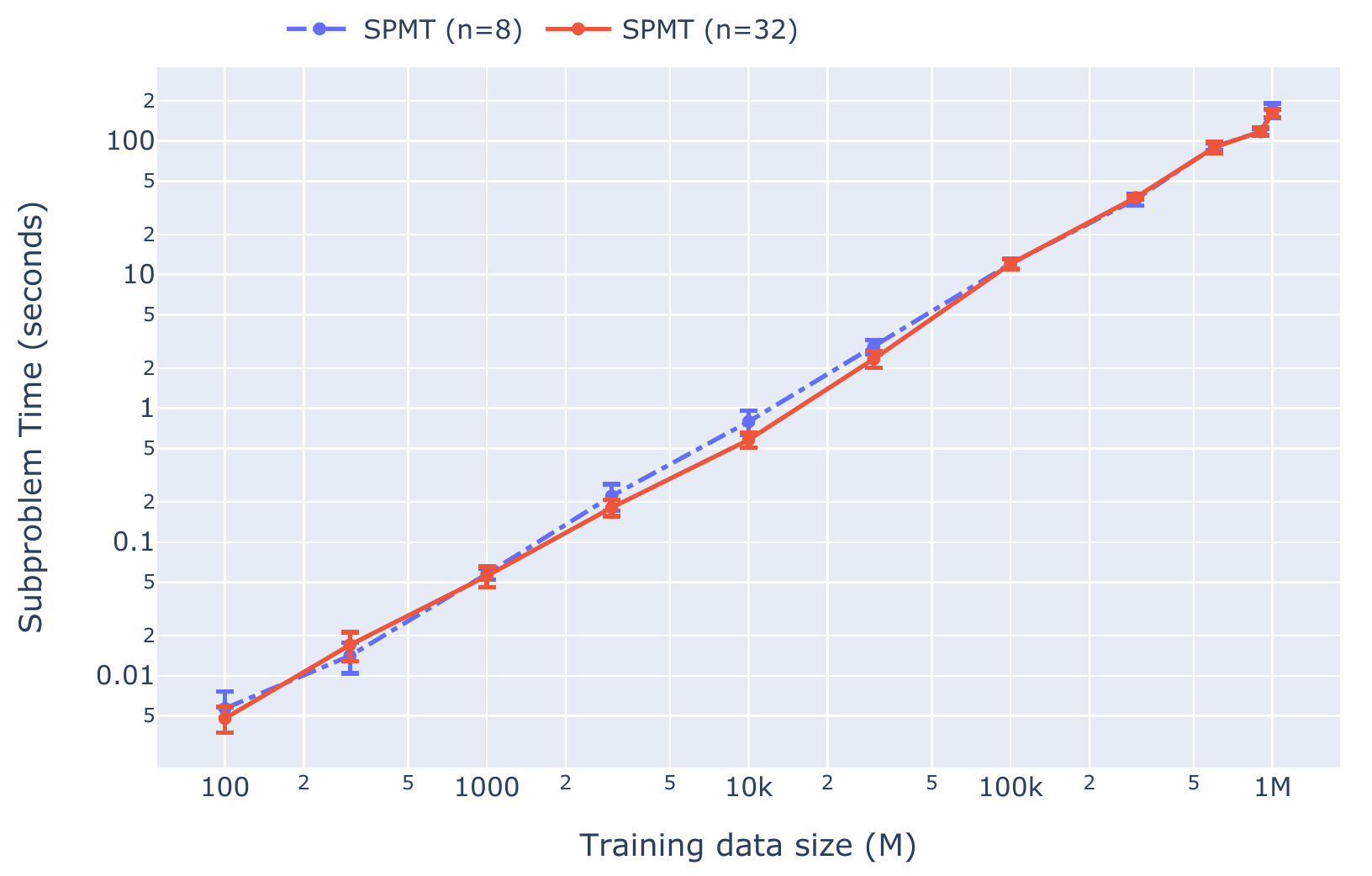}
  \caption{Average KSP duration}\label{rsub_time}
  \end{subfigure}
  \caption{Runtime results for synthetic Dataset (6)}\label{runtime_with_datasize}
\end{figure}

The log-log plots in Fig~\ref{runtime_with_datasize} show the average time taken to solve RMP and KSP, normalized by the number of CG iterations taken to converge. These metrics represent the computational effort required per CG iteration, independent of the data-dependent ``difficulty'' of the problem, which determines whether fewer or more CG iterations are required to converge. In our experiments, all instances took less than 20 CG iterations to converge, with 11.19 and 1.33 seconds as the mean and standard deviation. 

We do not observe significant change in these two metrics when we increase the cardinality of the rule set $n$. In fact, allowing more rules tends to reduce the runtime, as less work is required to ``finetune'' the rule set. 
As the samples increase, we observe a super-linear rise in the runtime per CG iteration that approximately follows a $M\log M$ trend. The average KSP runtime for the million sample data set was about 166 seconds while the RMPs were usually solved within 45 seconds by the CPLEX Barrier solver, enabling the CG algorithm to converge to a near optimal solution in about an hour.

\begin{figure}
  \begin{subfigure}{.5\textwidth}
  \includegraphics[width=1.1\linewidth]{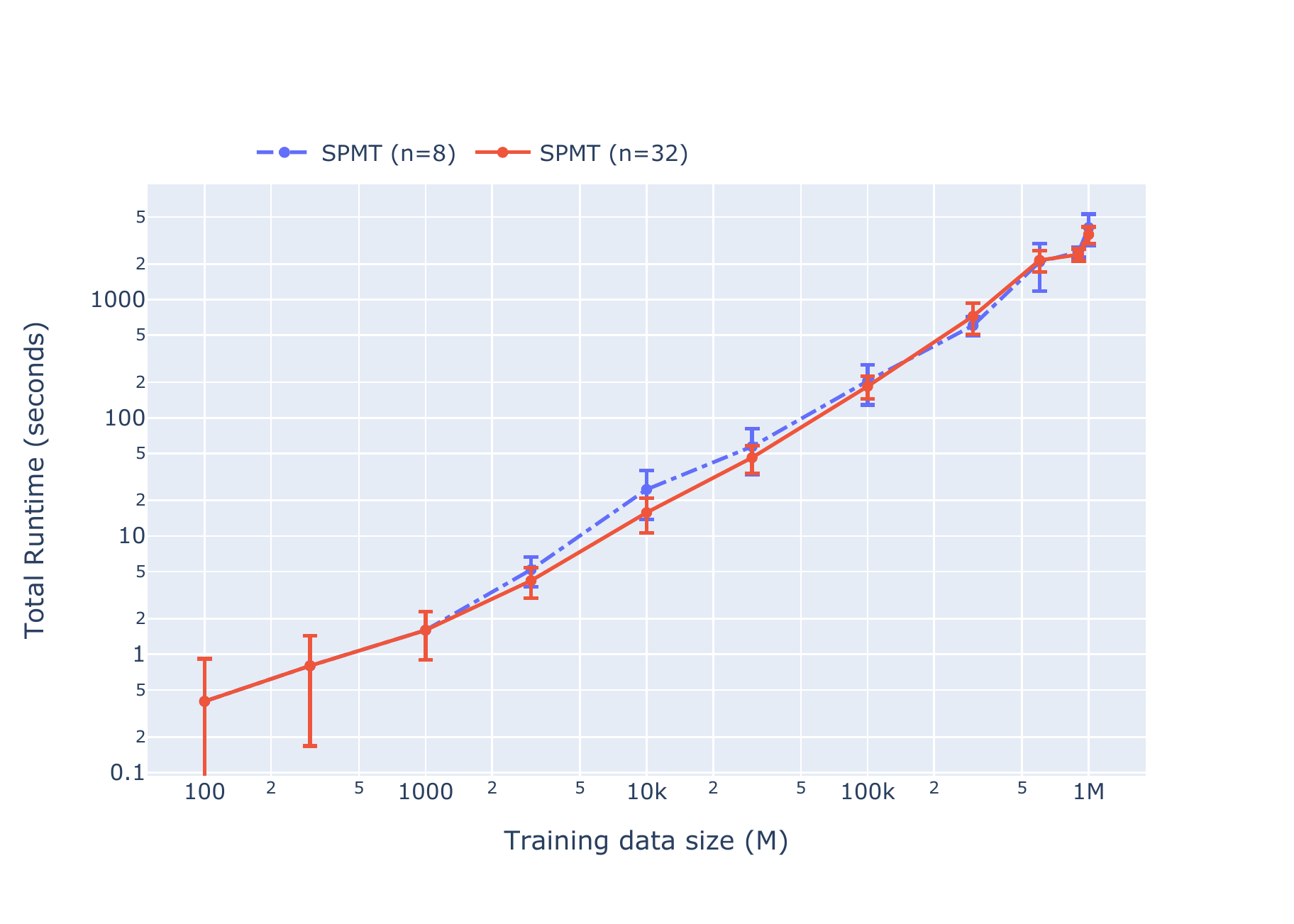}
  \caption{Total runtime}
  \end{subfigure}
\begin{subfigure}{.5\textwidth}
  \includegraphics[width=1.1\linewidth]{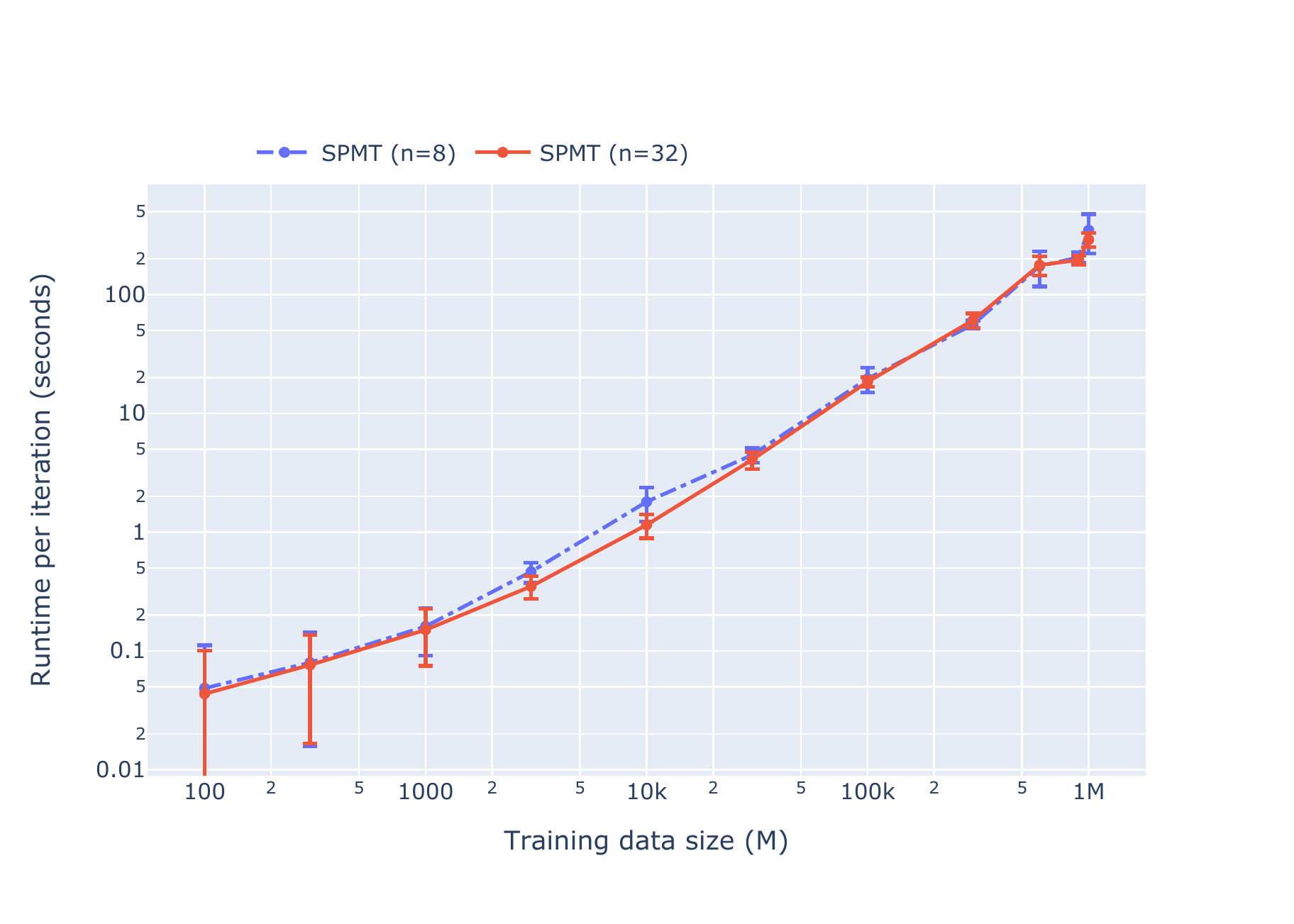}
  \caption{Average runtime per CG iteration}
  \end{subfigure}
\begin{subfigure}{.5\textwidth}
  \includegraphics[width=1.1\linewidth]{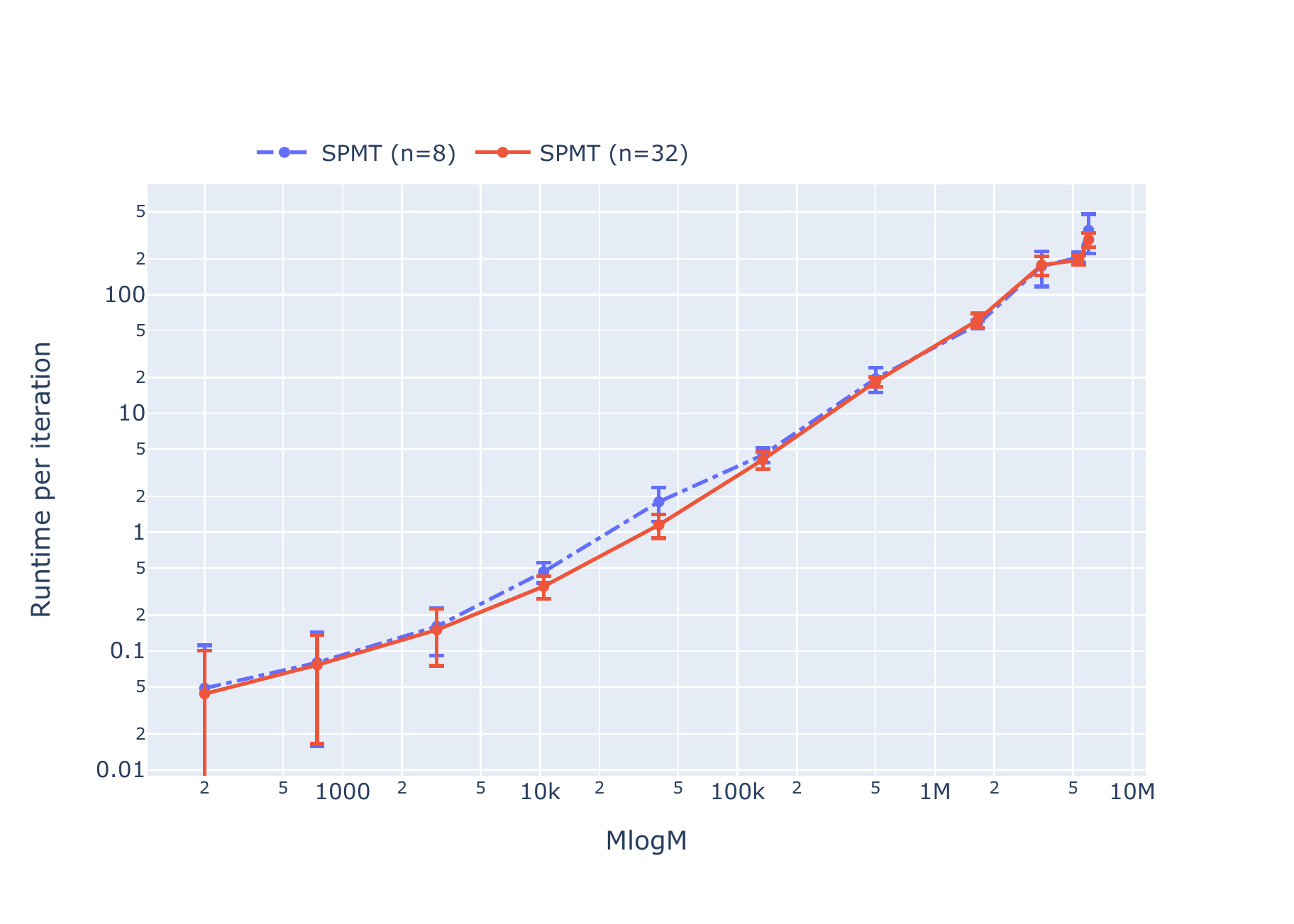}
  \caption{Average runtime against $MlogM$}
  \end{subfigure}
  \caption{Runtime experiment for Dataset (6)}
  \label{fig_runtime}
\end{figure}

Figure \ref{fig_runtime} provides additional information on the runtime experiments on Dataset (6). More specifically, total runtime records the amount of time (in seconds) to solve an instance, i.e., I/O operations, column generation (CG) iterations until convergence (including time taken to solve RMP and the subproblem), and the final  Master-MIP to obtain an integral solution. Average runtime is the total runtime divided by the number of CG iterations.  Note that
%
\shi{optimizing the model exactly} is an NP-Hard problem, which can be time-consuming for large, challenging problem instances. For simplicity and replicability, we solve \shi{our final Master-}MIP directly using a standard optimization library (\citealp{cplex2020}) to distill a (near) optimal subset of prescriptive decision rules.

The key takeaway from the runtime study is that our SPMT approach is more scalable \shi{in practical settings} than prior MIP approaches such as \cite{bertsimas2019optimal,amram2020optimal}, since the latter's MIP model  requires $O(2^kM)$ binary variables. For example, in the experiment which allowed a maximum of 32 rules for a dataset with $10^6$ training samples, prior MIP models would require more than $32\times10^6$ binary variables \shi{to be added upfront to obtain any answer}. In comparison, the Master-MIP in SPMT required less than 3000 binary variables \shi {to obtain near-optimal solutions}.

\begin{figure}
 \centering
  \includegraphics[width=0.4\linewidth]{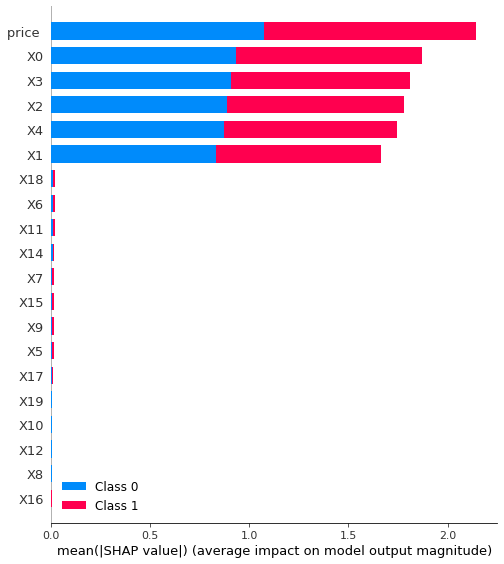}
  \caption{Feature importance based on the teacher model trained on Dataset (2)}
  \label{fig_shap}
\end{figure}

To demonstrate the efficacy of  column generation  in finding high quality solutions for large instances that have several million possible rules, we analyzed ten random instances of Dataset (2) containing $10^6$ data samples each for the multiway-split tree with $n=8$ and $32$ respectively. Dataset (2) was chosen for this study because of its relatively 
\shi{larger} feature space. Prior to constructing the tree, we analyzed the feature importance of the teacher model using SHAP package (\citealp{lundberg2017unified}), which correctly identifies five features along with price as the most predictive features on the outcome, which are used to construct the feature graph (see Fig~\ref{fig_shap}). This highlights a benefit of decoupling prediction from the downstream policy optimization task, as the prescriptive student is able to utilize the teacher's knowledge to efficiently filter out the uninformative features and construct a more compact feature graph. 

\begin{figure}
 \centering
  \includegraphics[width=0.65\linewidth]{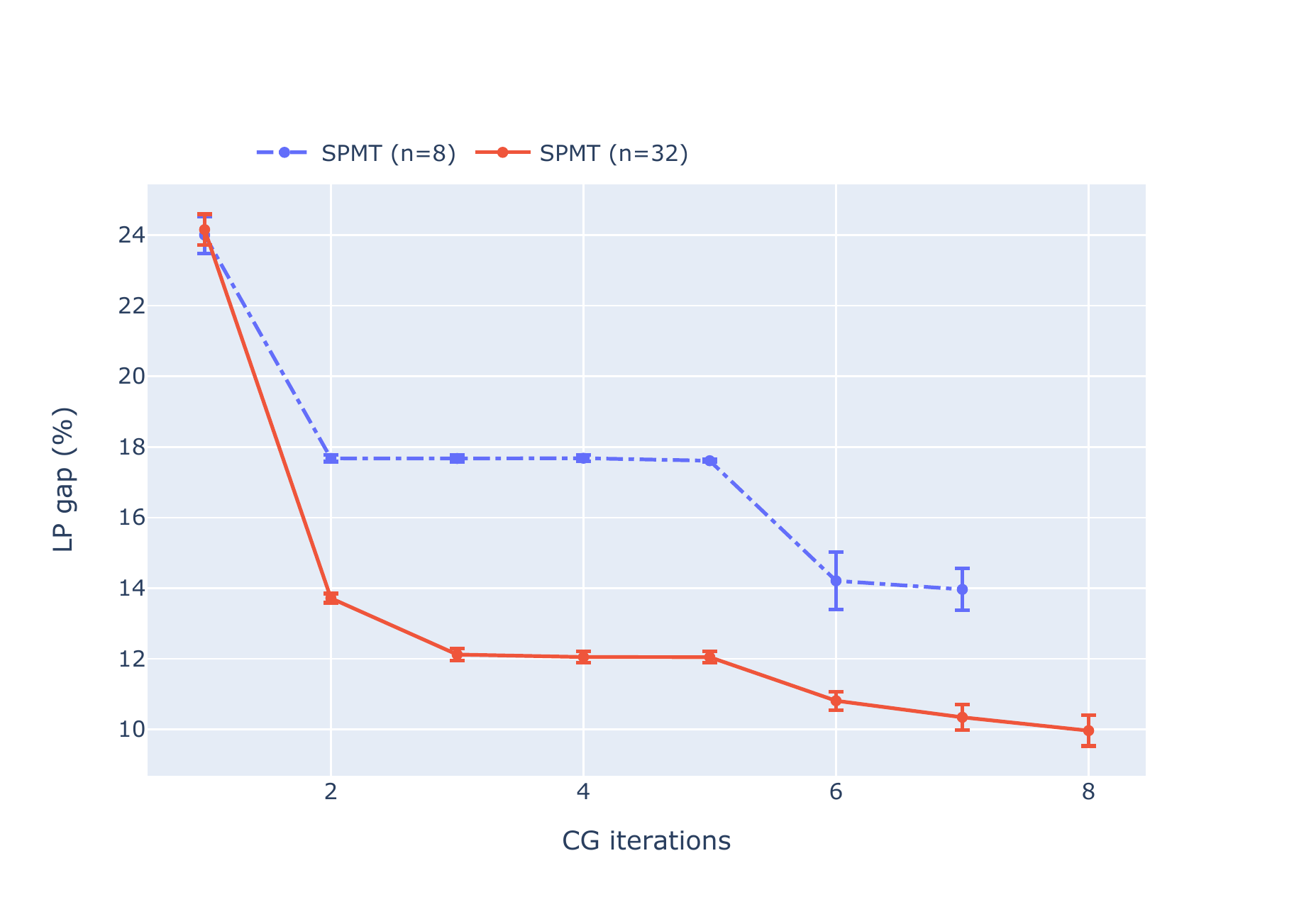}
  \caption{Progression of column generation for Dataset (2)}
  \label{fig_CG_progression}
\end{figure}

Each of the five numerical features were discretized into 12 bins ($\kappa = 4$) using the cumulative binning method discussed earlier. The number of possible paths for this setting in combination with 9 price levels is $13^5\times9$, which results in more than 3.3 million possible rules. We compare the relative gap between between the RMP objective function and the teacher model predicted outcome (obtained by selecting the revenue-maximizing price for each sample). Note that the predicted teacher outcome  acts as \shi{a global} upper bound (UB), which may not be achievable for a fixed $n$ where $n$ is much smaller than $M$. For this experiment,
we set $K = 500$ and a maximum time limit of one hour for the CG iterations to reduce the relative gap between the \shi {achieved} RMP objective function and the UB to 10\%. If the CG does not converge within the time limit, we terminate the column generation and solve the Master-MIP. 
Meanwhile, no time limit is imposed on solving the resultant Master-MIP. The achieved RMP objective function value after each CG iteration is depicted in Fig~\ref{fig_CG_progression}. 

Table~\ref{tab_data2} reports on the average values of achieved MIP and LP (RMP) objective function values expressed in terms of the percentage gap from 
\shi{UB}, followed by the final number of $z$-variables (RMP columns), and the average run time (in minutes) for solving the RMP, the subproblem and the Master-MIP. The standard deviations over ten instances are shown in parentheses. This table shows that with $n=32$ rules, the solution is within 10.3\% of UB, whereas the average final gap that permits no more than 8 rules is 17\%. We observe small differences between this final LP and MIP gaps \shi {from UB} that is less than 0.17\% and 0.61\%  with $n=32$ and 8, respectively. \shi{Note that it is not necessary to solve the LP to provable optimality in practical settings where exact solutions are not required. In this case, the achieved LP objective is not a guaranteed bound on a global optimal solution but is a valid bound on the Restricted Master MIP objective.} This narrow MIP-LP gap highlights the benefit of converging to a relatively small pool of high quality rules from which an effective subset of final rules can be recovered via a standard MIP solver. 

A key finding in this experiment is that through column generation we required less than 3000 positive reduced cost paths in all instances to distill a near-optimal set with at most 32 rules. This translates to $\hat{N} < 0.0009N$, i.e. an RMP having less than 0.1\% of all possible binary variables, which represents a dramatic reduction in the  problem's effective complexity. Even fewer paths (less than 0.05\% of all possibilities) were required for $n=8$ rules. 
{In comparison, for $n=32$ rules, prior MIP approaches would require a model having more than 30 million binary variables to be formulated, which would be computationally intractable. }

In terms of runtime, Table~\ref{tab_data2} shows that a 10.3\% gap \shi {from UB} was achieved with $n=32$ in slightly more than an hour. The maximum time limit for CG iterations was reached in all instances for $n=32$. In contrast, the average run time for $n=8$ to attain a solution gap of 17\% \shi {from UB} was relatively lower (about 39 minutes) albeit with a higher standard deviation, as the CG iterations converged before the time limit in seven of the ten instances. However, as seen in Fig~\ref{fig_CG_progression}, the achieved LP gap \% with $n=32$ after just two CG iterations is better than the final LP gap achieved with $n=8$. A less restrictive cardinality constraint resulted in faster convergence to the same solution quality in our experiments.

\begin{table}
\centering
\caption{Worst-case optimality gap and runtime for Dataset (2)}
\label{tab_data2}
\begin{tabular}{c|cccc}
\toprule
 Rules ($n$) &  MIP gap \% &  LP gap \% &  RMP columns &  Average runtime (min) \\
\midrule
     8 &  17.033 (2.031) & 16.581 (1.827) &       811 (363) &              38.579 (16.972) \\
    32 &  10.281 (0.457) & 10.164 (0.452) &      2539 (251) &               65.351 (5.155) \\
\bottomrule
\end{tabular}
\end{table}


\subsection{Real datasets}\label{sect_real}\label{use_case_grocery}

We perform multiple experiments with real datasets to demonstrate the practical benefits of adding various intra-rule and inter-rule constraints to satisfy critical real-life requirements. Specifically, we are able to identify non-discriminatory and transparent prescriptive policies that also increase revenue, operational effectiveness, and proactively avoid inter-rule conflicts. We achieve this task in two steps. First, we execute the CG algorithm in an ``unconstrained'' mode, where we do not consider additional restrictions beyond the rule-cardinality constraint. Subsequently, in one or more ``constrained'' modes, we impose additional side constraints and report on the performance improvements with respect to the different goals.

\subsubsection{Grocery pricing case study}

One potential issue with the unconstrained pricing policy is its significantly high demand predicted by the teacher model. 
Compared to the median historical price of $\$2.99$ for 1 pound of strawberries, 
the bulk of the optimal policies suggest a large price cut to \$2 as shown in Fig 1b in the main paper. The key insight is that the demand for strawberries is elastic across many segments, i.e.,  a decrease in price leads to a significant increase in demand and its revenue. In fact, the predicted demand at the optimal prescribed policy is 99.94\% higher than the historical value. Such demand surge may be neither  realistic nor desirable
due to supply constraint or other limits imposed by the perishable nature of the product.  Such observations motivate us to  investigate this use case by imposing realistic operational constraints.

We derive two features to facilitate store-cluster pricing and loyalty pricing. More specifically, we create a feature called ``Cluster'' as we group 311  stores in the dataset into 5 clusters, i.e., ``extremely low'', ``low'', ``medium'', ``high'' and ``mega'' sales segments,  based on their past sales by using the following cutoff: $[0, 2, 15, 30, 100, 300]$. That is, a store with total sales within [15, 30) will be assigned to ``medium'' sales cluster. Note that in this dataset, 199 stores fall into the ``extremely low'' sales cluster. Meanwhile, we also create a binary indicator ``Loyalty'' to indicate whether the customer is a loyalty card holder. To achieve this, we analyze the original Dunnhumby dataset\footnote{\url{https://www.dunnhumby.com/careers/engineering/sourcefiles}}. For each customer, Loyalty = 1 if one of her historical transactions (not restricted to strawberry purchases) contains a non-zero ``retail\_disc'' value which indicates a retail discount
applied due to a loyalty card. Based on this definition, in the processed dataset\footnote{\url{https://docs.interpretable.ai/stable/examples/grocery_pricing/}}  on strawberry purchases, 465 out of 547 customers  are identified as loyalty card holders.

With store-cluster pricing, we impose a constraint for each cluster which requires the predicted demand under a new pricing policy does not deviate too much from the historical average.  With 64 rules, SPMT produces policies that result in 47.49\% and 61.72\% increase in revenue over the historical policy when the allowable change in predicted demand is set to 25\% and 50\% of the historical value for each cluster. Unsurprisingly, more stringent capacity constraint results in lower revenue. 

With loyalty price, we use StoreID and Loyalty indicator as features for the tree, where all shoppers at a store receive the same price except loyalty-card shoppers, who receive a price that is same or lower. 
%
After applying these constraints, SPMT recommends a policy which leads to 65.4\% improvement in revenue over the historical baseline. 


\subsubsection{Housing upgrade use case}


This use case is different from other studies in the paper because the teacher model addresses a regression problem, i.e.,  predict the sales price of a house given its attributes, as opposed to being a classifier. We assume that the houses in the evaluation data are yet to be constructed and their grades can be altered to change their predicted sales prices. The predicted gain is given by the sum of the teacher predicted sales prices. 
Numerical features such as living area and age were aggregated into bins and the CG algorithm was run on the resulting data set. The output decision rules were then evaluated on the properties listed in the test data.

We first analyze the unconstrained setting. The unconstrained policy yields a single decision rule to upgrade to the highest level, which essentially increases the predicted sales prices of all houses to their highest value and results in a 124\% gain in the test data. However, such a rule destroys the price diversity of the houses in a neighborhood and ignores that the incremental investment required to upgrade all houses. The other two scenarios address these concerns using constraints.

A key feature of our approach is the ability to manage constraints to ensure fair and operationally effective prescriptive rules. Using a constrained housing scenario as an example, we explain how inter-rule or global constraints are satisfied. Suppose we impose a ``zip-4'' pricing constraint at the neighborhood level that ensure that the average predicted sales price of houses having the same first 4 digits of the zip code cannot change by more than 10\%. There are a total of 19 such localities in the dataset and this restriction is applied to each of them as discussed below.

Such a constraint can be  added to the RMP, resulted in a slightly different reduced cost to guide column generation. Denote $d_{il}=1$ if house $i$ is in zipcode $l$, 0 otherwise, $p_l$ as the historical sale price for  zipcode $l$. Recall $g_{i,q_j}$ is the predicted counterfactual outcome (sale price) associated with rule $j$ for sample (house) $i$. The required constraint can be formulated as,  
\begin{equation}\label{housing_constraint}
    0.9p_l\leq \frac{1}{\sum_i d_{il}}\sum_{j}\sum_i d_{il}g_{i,q_j}z_j \leq 1.1p_l, \quad\forall i, l \tag{A.1}
\end{equation}

Eq~\ref{housing_constraint} is simply a linear constraint as $d_{il}$ and $g_{i,q_j}$ are parameters to the optimization problem.   We can rewrite Eq~\ref{housing_constraint} as $$\frac{1}{\sum_i d_{il}}\sum_{j}\sum_i d_{il}g_{i,q_j}z_j \leq 1.1p_l, $$ and $$-\frac{1}{\sum_i d_{il}}\sum_{j}\sum_i d_{il}g_{i,q_j}z_j \leq -0.9p_l.$$ Let $\rho^t_{l+}$ and $\rho^t_{l-}$ denote the non-negative dual variables associated with these two sets of constraints respectively. Then the reduced cost can be written as
$$rc_j^t =  r_j^t -\left(\sum_{i}a_{ij}\lambda_i^t + \mu^t+ \sum_{i }\sum_{l}\frac{1}{\sum_i d_{il}}d_{il}g_{i,q_j}(\rho^t_{l+}-\rho^t_{l-})  \right).$$

More generally, inter-rule constraints which span across several paths/rules (e.g., the cardinality constraint as Eq (4) in the main paper and Eq~\ref{housing_constraint}) enter the RMP formulation. They influence the subproblem through the resultant reduced costs to ensure that paths that are more likely to be feasible are added to the MIP.  Inter-rule restrictions can be specified as convex or linear (continuous or discrete) inequalities. Constraint satisfaction is achieved in the final Master-MIP when we reimpose the binary restrictions on $z$.

Imposing this constraint limits the achievable revenue gain while generating a grading policy that preserves the existing price characteristics of each locality. The decision rules no longer upgrade all houses but contain a mixture of upgrades, downgrades, or no changes based on the living area, age, and number of baths and yields a predicted revenue gain of 9.74\%.

In the second constrained scenario which bears some resemblance to a portfolio rebalancing case, we assume the cost of an upgrade is taken as 70\% of the predicted sales price increase along with a similar cost saving by downgrading. No cost is incurred for leaving the grade unchanged. The goal is to generate no more than $8$ decision rules for grades such that the net cost increase from all the grade changes globally is close to zero.  This requirement is practically enforced by limiting the net cost change to less than 1\% of the baseline cost (taken as 70\% of the predicted sales prices at the original grades). The resulting rules when scored on the test data set yields a predicted revenue gain of 7.82\%. The rules essentially rationalize the grades with respect to the predicted return on investment and recommend diverse upgrades and downgrades primarily based on the aggregate living area of the house.




\subsubsection{Airline pricing use case}

The dataset has more than $1.3\times 10^6$ samples of buy and no-buy data spanning several markets containing more than 10 different multi-level categorical features, and 15 price points on a price grid to choose from. \shi{An example of an acceptable rule may look like the following: ``if the flight is between Airports A and B, and the ticket was purchased less then 7 days before departure, and the departure time is on a weekday before 11am, then a price of \$400 is prescribed".} The column generation algorithm was run in two modes. In the first mode, we ignored capacity constraints and achieved a overly high predicted revenue gain that was within 12\% of the upper bound (UB) obtained by identifying an optimal personalized price rule for every transaction. Note that these personalized prescriptions potentially charge different prices for different customers who buy a ticket at the same time, which violates the user-acceptance requirements. The unconstrained policy satisfies this requirement but creates other operational problems.  

Analysis revealed that the unconstrained gain can only be achieved by selling more than twice the number of tickets being sold currently, which is dramatically higher than the available capacity. Airline travelers are price sensitive, and conversion rates are quite sensitive to price changes (e.g., a 10\% price cut can result in a 25\% increase in conversions for some customer segments). If such rules are deployed it will allow leisure travelers who book early to grab premium seats at a steep bargain, which in turn can result in early stock-outs. Premium seats will become unavailable to high value travelers who book later and drive them toward competitor airlines in that market. By adding carefully calibrated market-specific constraints on the teacher-model predicted conversions at different prices, we are able to \shi{re-optimize} the rules satisfactorily and still achieve a predicted gain of more than 25\% that is more likely to be realized in live testing compared to the unconstrained policy. 

In terms of runtime performance for this challenging constrained setting, the CG iterations converged in about 45 minutes, while the final MIP with the complicating market-specific capacity constraints was solved using the CPLEX solver in a few hours. As the available capacities can change periodically, the model can be warm-started from its most recent solution to quickly re-optimize the decision rules on demand.

\bibliographystyle{plainnat}
\bibliography{ref}

